\documentclass[a4paper,10pt]{amsart}
\usepackage{fourier}
\usepackage[margin=2.75cm]{geometry}
\usepackage{mathrsfs}
\usepackage{amssymb}
\usepackage{color}
\usepackage[colorlinks=true,linkcolor=blue,citecolor=blue]{hyperref}
\usepackage{xy}
\usepackage{xypic}
\usepackage{pgf,tikz}

  \usetikzlibrary{arrows}
  \usetikzlibrary{decorations.markings}
  \usetikzlibrary{shapes}
  \usetikzlibrary{matrix}
\usepackage{wasysym}

\newcommand{\C}{\mathcal{C}}
\newcommand{\B}{\mathcal{B}}
\newcommand{\A}{\mathcal{A}}
\newcommand{\X}{\mathscr{X}}
\newcommand{\E}{\mathscr{E}}
\newcommand{\Y}{\mathscr{Y}}

\renewcommand{\diagram}[3]{\matrix (#1) [matrix of math nodes,row
  sep={#2},column sep={#3},text height=1.5ex,text depth=0.25ex]}

\DeclareMathOperator{\Hom}{Hom}

\numberwithin{equation}{section}
\def\P{\mathscr{P}}
\def\I{\mathscr{I}}
\theoremstyle{plain}
\newtheorem{thm}{subsection}[section]
\newtheorem{theorem}[thm]{Theorem}
\newtheorem{proposition}[thm]{Proposition}
\newtheorem{lemma}[thm]{Lemma}
\newtheorem{corollary}[thm]{Corollary}
\newtheorem{example}[thm]{Example}
\theoremstyle{definition}
\newtheorem{definition}[thm]{Definition}
\newtheorem{remark}[thm]{Remark}

\title{One-sided $n$-suspended categories}
\author{Jing He, Yonggang Hu and Panyue Zhou }

\address{Department of Mathematical Sciences, Tsinghua University, 100084 Beijing,  P. R. China}
\email{huyonggang@emails.bjut.edu.cn}

\address{College of Mathematics, Hunan Institute of Science and
Technology, 414006, Yueyang, Hunan,
	P. R. China.}
\email{panyuezhou@163.com}

\begin{document}

% Abstract
\begin{abstract}
Let $n$ be
an integer greater or equal than $3$. We give a simultaneous generalization of $(n-2)$-exact categories and $n$-angulated
categories, and we call it  one-sided $n$-suspended categories. One-sided $n$-angulated categories
are also examples of  one-sided $n$-suspended categories. We
provide a general framework for passing from one-sided $n$-suspended categories to one-sided $n$-angulated categories. Besides, we give a method to construct $n$-angulated quotient categories from Frobenius $n$-prile categories.  These results generalize their works by Jasso for $n$-exact categories, Lin for $(n+2)$-angulated categories and  Li for one-sided suspended categories.
\end{abstract}
%\thanks{$^{\ast}$~Corresponding author.  }
\thanks{Yonggang Hu was supported
	by the National Natural Science Foundation of China (Grant Nos. 11971255 and 12071120). Panyue Zhou
was supported by the National Natural Science Foundation of China (Grant No. 11901190) and
by the Scientific Research Fund of Hunan Provincial Education Department (Grant No. 19B239). }


\maketitle

\section{Introduction}
Geiss, Keller and Oppermann \cite{[GKO]} introduced the notion of $n$-angulated category, which is
a generalization of triangulated category.
Certain $(n-2)$-cluster tilting subcategories of triangulated categories
are examples of $n$-angulated categories \cite{[GKO]}.
 Recently, Jasso \cite{[J]} introduced the notions of  $(n-2)$-exact
categories and algebraic $n$-angulated categories, where
algebraic $n$-angulated categories are by definition equivalent to stable categories of  Frobenius $(n-2)$-exact categories, which are natural generalizations of algebraic triangulated categories.
Lin \cite{[L2]} defined mutation pairs in $n$-angulated categories and proved that given such
a mutation pair, the corresponding quotient category is an $n$-angulated category.
\vspace{1mm}

The main aim of this paper is to provide a general framework to understand the $n$-angulated
quotient categories appeared in \cite{[J],[L2]}. Inspired by the notions of one-sided $n$-angulated category introduced by Lin \cite{[L1]} and one-sided suspended category by Li \cite{[Li]}, we define one-sided $n$-suspended categories, which contain $(n-2)$-exact categories, $n$-angulated categories and one-sided $n$-angulated categories as examples.
\vspace{1mm}

Roughly speaking, a  right $n$-suspended
category consists of an additive category $\A$ endowed with an additive endofunctor $\Sigma$, two additive
subcategories $\X$ and $\C$ of $\A$ together with a class $R_n(\C,\Sigma)$ of $n$-$\Sigma$-sequences, which satisfy similar axioms of a right $n$-angulated category \cite{[L1]}. Right $n$-angulated categories and $(n-2)$-exact
categories are examples of right $n$-suspended categories. The precise definition of a
right $n$-suspended category can be found in Definition \ref{d2}. A  left $n$-suspended category is defined dually. Our first main result is the following:

\begin{theorem}(see Theorem \ref{y9} for details) \emph{(1)} If $(\A,R_n(\C,\Sigma),\X)$ is a right $n$-suspended category, then the quotient category
$\C/\X$ is a right $n$-angulated category.
\vspace{1mm}

\emph{(2)} If $(\A,\Omega,L_n(\C),\X)$ is a  left $n$-suspended category, then the quotient category
$\C/\X$ is a left $n$-angulated category.

\end{theorem}

We use our first main result to give about the construction of $n$-angulated quotient categories. It covers many existed constructions of $n$-angulated quotient categories \cite[Theorem 5.11]{[J]} and \cite[Theorem 3.7]{[L2]}.  More precisely, we prove the following.

\begin{theorem}(see Theorem \ref{h11} for details)
Let $(\A,L_n(\C,\Omega),R_n(\C,\Sigma),\X)$ be a Frobenius $n$-prile category. Then
the quotient category $\C/\X$ is an $n$-angulated category.
\end{theorem}

This paper is organized as follows. In Section 2, we recall the definitions of $(n-2)$-exact category and right $n$-angulated category.  In Section 3,  we define the notion of a  one-sided $n$-suspended category and give some constructions of  one-sided $n$-suspended category by  $(n-2)$-exact categories and one-sided $n$-angulated categories. In Section 4, we provide a general framework for passing from a one-sided $n$-suspended category to an one-sided $n$-angulated category.
 In Section 5, we introduce the notion
of a Frobenius $n$-prile category and give some examples, then prove our second main result.
\vspace{1mm}

\textbf{Conventions}. Throughout this paper, $n$ is an integer greater than or equal to three.
Let $\A$ be an additive category, when we say that $\X$
is a subcategory of $\C$, we always mean that $\X$ is full and is closed under isomorphisms, direct sums, and
direct summands.

\section{Partial one-sided $n$-angulated categories}
\subsection{$(n-2)$-exact categories}
In this subsection, we recall some definitions of an $(n-2)$-exact category from \cite{[J]}.
\vspace{1mm}

Let $\A$ be an additive category and $f\colon A\rightarrow B$ a morphism in $\A$. A \emph{weak cokernel} of $f$ is a morphism
$g\colon B\rightarrow C$ such that for any $X\in\A$ the sequence of abelian groups
$$\A(C,X)\xrightarrow{g^\ast}\A(B,X)\xrightarrow{f^\ast}\A(A,X)$$
is exact. Equivalently, $g$ is a weak cokernel of $f$ if $gf=0$ and for each morphism
$h\colon B\rightarrow X$ such that $hf=0$ there exists a (not necessarily unique) morphism
$p\colon C\rightarrow X$ such that $h=pg$. These properties are subsumed in the following
commutative diagram
$$\xymatrix{A\ar[r]^{f}\ar[dr]_{0} & B\ar[r]^{g}\ar[d]^{\forall h}&C\ar@{-->}[dl]^{\exists p}  \\
& X &}$$
Clearly, a weak cokernel $g$ of $f$ is a cokernel of $f$ if and only if $g$ is an epimorphism.
The concept of a \emph{weak kernel} is defined dually.
\vspace{1mm}

Let $f_1\colon A_1\rightarrow A_2$ be a morphism in
$\A$. An $(n-2)$-\emph{cokernel} of $f_1$ is a sequence
$$(f_2,f_3,\cdots,f_n)\colon A_2\xrightarrow{f_2}A_3\xrightarrow{f_3}\cdots\xrightarrow
{f_{n-1}}A_n$$
such that for any $B\in\A$ the induced sequence of abelian groups
$$\xymatrix@C=1cm{0\xrightarrow{}\A(A_{n},B)\xrightarrow{(f_{n-1})_\ast} \A(A_{n-1},B)\xrightarrow{(f_{n-2})_\ast}\cdots\xrightarrow{(f_1)_\ast}
\A(A_{1},B)}$$
is exact. Equivalently, the sequence $(f_2,f_3,\cdots,f_{n-1})$ is an $(n-2)$-cokernel of $f_1$ if for any
$2\leq i\leq n-2$ the morphism $f_i$ is a weak cokernel of $f_{i-1}$, and $f_{n-1}$ is moreover a
cokernel of $f_{n-2}$. In this case, we say the sequence
\begin{equation}\label{t0}
\begin{array}{l}
A_1\xrightarrow{f_1}A_2\xrightarrow{f_2}A_3\xrightarrow{f_3}\cdots\xrightarrow
{f_{n-1}}A_n
\end{array}
\end{equation}
is a \emph{right $(n-2)$-exact}.
We can define $(n-2)$-\emph{kernel} and left $(n-2)$-\emph{exact sequence} dually.
The sequence (\ref{t0}) is called $(n-2)$-\emph{exact} if it is both right $(n-2)$-exact and left $(n-2)$-exact.

Let $\A$ be an additive category, $A_1\xrightarrow{f_1} A_2\xrightarrow{f_2}\cdots\xrightarrow{f_{n-2}}A_{n-1}$ a complex and $h_1\colon A_1\rightarrow B_1$ a morphism in $\A$. An $(n-2)$-\emph{pushout} diagram of $(f_{1}, f_{2}, \cdots, f_{n-2})$ along $h_1$ is a morphism of complexes
$$\xymatrix{
&A_1\ar[r]^{f_1}\ar[d]^{h_1} & A_2\ar[r]^{f_2}\ar[d]^{h_2} & \cdots \ar[r]^{f_{n-3}} & A_{n-2}\ar[r]^{f_{n-2}}\ar[d]^{h_{n-2}}  & A_{n-1} \ar[d]^{h_{n-1}}\\
&B_1\ar[r]^{g_1} & B_2\ar[r]^{g_2}  & \cdots \ar[r]^{g_{n-3}} & B_{n-2}\ar[r]^{g_{n-2}} &  B_{n-1}\\
}$$
such that the mapping cone
$$A_1\xrightarrow{\left[
                                            \begin{smallmatrix}
                                              -f_1 \\
                                              h_1 \\
                                            \end{smallmatrix}
                                          \right]}
 A_2\oplus B_1\xrightarrow{\left[
                                            \begin{smallmatrix}
                                              -f_{2} & 0 \\
                                              h_{2} & g_1
                                            \end{smallmatrix}
                                          \right]}
A_3\oplus B_2\xrightarrow{\left[
                                            \begin{smallmatrix}
                                              -f_{3} & 0 \\
                                              h_{3} & g_2
                                            \end{smallmatrix}
                                          \right]}
\cdots\xrightarrow{\left[
                                            \begin{smallmatrix}
                                              -f_{n-2} & 0 \\
                                              h_{n-2} & g_{n-3} \\
                                            \end{smallmatrix}
                                          \right]\ \ \ }
A_{n-1}\oplus B_{n-2}\xrightarrow{\left[\begin{smallmatrix}
                                             h_{n-1} & g_{n-2} \\
                                            \end{smallmatrix}
                                          \right]}
B_{n-1}$$
is a right $(n-2)$-exact.
\medskip

Let $\A$ be an additive category. If $\X$ is a class of $(n-2)$-exact
sequences in $\A$, then we call its members $\X$-\emph{admissible
  $(n-2)$-exact sequences}.  Furthermore, if
$$A^{\bullet}:=(A_1\xrightarrow{f_1}A_2\xrightarrow{f_2}A_3\xrightarrow{f_3}\cdots\xrightarrow
{f_{n-1}}A_n)$$
is an $\X$-admissible $(n-2)$-exact sequence, we say that $f_1$ is an
\emph{$\X$-admissible monomorphism} and that $f_{n-1}$ is an
\emph{$\X$-admissible epimorphism}.  When the class $\X$ is clear from
the context, we write ``admissible'' instead of ``$\X$-admissible''.

\begin{definition} {\cite[Definition 4.2]{[J]}}
We say that a morphism $f^{\bullet}\colon A^\bullet\rightarrow B^\bullet$ of $(n-2)$-exact sequences in an additive category $\A$ is a \emph{weak isomorphism} if $f^i$ and $f^{i+1}$ are isomorphisms for some $i\in\{1,2,\cdots, n\}$ with $n+1:=1$.  An \emph{$(n-2)$-exact structure} on $\A$ is a class $\X$ of $(n-2)$-exact sequences in $\A$, closed under weak
isomorphisms of $(n-2)$-exact sequences, and which satisfies the following axioms:

\begin{itemize}
\item[(E0)] The sequence $0\rightarrow  0\rightarrow \cdots \rightarrow 0\rightarrow 0$ is an admissible $(n-2)$-exact sequence.
\item[(E1)] The class of admissible monomorphisms is closed under composition.
\item[(E1)$^{\textrm{op}}$] The class of admissible epimorphisms is closed under composition.
\item[(E2)] For any admissible $(n-2)$-exact sequence $A^{\bullet}$ and any morphism $\varphi\colon A_1\rightarrow B_1$, there exists an $(n-2)$-pushout diagram of $(f_1,f_2,\cdots,f_{n-2})$
along $\varphi$ such that $g_1$ is an admissible monomorphism, that is, there exists the
following commutative diagram:
$$\xymatrix{
A_1\ar[r]^{f_1}\ar[d]^{\varphi} & A_2\ar[r]^{f_2}\ar@{-->}[d]^{h_2} & A_3\ar[r]^{f_3}\ar@{-->}[d]^{h_3} & \cdots \ar[r]^{f_{n-2}} &A_{n-1}\ar@{-->}[d]^{h_{n-1}} \ar[r]^{f_{n-1}}&A_n \\
B_1\ar[r]^{g_1} & B_2\ar@{-->}[r]^{g_2} & B_3\ar@{-->}[r]^{g_3}  & \cdots \ar@{-->}[r]^{g_{n-2}} & B_{n-1}
}$$
\item[(E2)$^{\textrm{op}}$] For any admissible $(n-2)$-exact sequence $D^{\bullet}$ and any morphism $\psi\colon C_{n}\rightarrow D_{n}$, there exists an $(n-2)$-pullback diagram of $(d_2,d_3,\cdots,d_{n-1})$ along $\psi$ such
that $c_{n-1}$ is an admissible epimorphism, that is, there exists the following commutative diagram:
$$\xymatrix{
& C_2\ar[r]^{c_2}\ar@{-->}[d]_{e_2} & C_3\ar@{-->}[r]^{c_3}\ar@{-->}[d]^{e_3} & \cdots \ar[r]^{c_{n-2}} & C_{n-1}\ar[r]^{c_{n-1}}\ar@{-->}[d]^{e_{n-1}}&C_n\ar[d]^{\psi}\\
D_1\ar[r]^{d_1} & D_2\ar[r]^{d_2} & D_3\ar[r]^{d_3}  & \cdots \ar[r]^{d_{n-2}} & D_{n-1}\ar[r]^{d_{n-1}}&D_{n}}$$
\end{itemize}
An $(n-2)$-\emph{exact category} is a pair $(\A,\X)$ where $\A$ is an additive category and $\X$ is an
$(n-2)$-exact structure on $\A$. If the class $\X$ is clear from the context, we identify $\A$ with
the pair $(\A,\X)$.
\end{definition}

Let $(\A,\X)$ be an $(n-2)$-exact category. An object $I\in\A$ is called
\emph{injective} if for every admissible monomorphism $f\colon M\to N$, the sequence of abelian groups
$$\xymatrix{\A(N,I)\ar[r]^{f^\ast}&{\A(M,I)}}\longrightarrow 0$$
is exact. By \cite{[J]},  we say that $(\A,\X)$ \emph{has enough injectives}
if for every object $M\in\A$ there exists injective objects
$I_2,\dots,I_{n-1}$ and an admissible $(n-2)$-exact sequence
$$M\rightarrow I_2\longrightarrow \cdots \longrightarrow I_{n-1}\rightarrow N.$$
The notion of \emph{having enough projectives} is defined dually.
\vspace{1mm}

We say that an $(n-2)$-exact category $(\A,\X)$ is \emph{Frobenius} if
  it has enough injectives, enough projectives, and if
  injective and projective objects coincide.

\subsection{Right $n$-angulated categories}
Let $\mathcal{A}$ be an additive category with an endofunctor $\Sigma:\mathcal{A}\rightarrow\mathcal{A}$. An $n$-$\Sigma$-$sequence$ in $\mathcal{A}$ is a sequence of morphisms
$$A_1\xrightarrow{f_1}A_2\xrightarrow{f_2}A_3\xrightarrow{f_3}\cdots\xrightarrow{f_{n-1}}A_n\xrightarrow{f_n}\Sigma A_1.$$
Its {\em left rotation} is the $n$-$\Sigma$-sequence
$$A_2\xrightarrow{f_2}A_3\xrightarrow{f_3}A_4\xrightarrow{f_4}\cdots\xrightarrow{f_{n-1}}A_n\xrightarrow{f_n}\Sigma A_1\xrightarrow{(-1)^n\Sigma f_1}\Sigma A_2.$$  %We can define $right\ rotation\ of\ an$ $n$-$\Sigma$-$sequence$ similarly.
A \emph{morphism} of $n$-$\Sigma$-sequences is  a sequence of morphisms $\varphi=(\varphi_1,\varphi_2,\varphi_3,\cdots,\varphi_n)$ such that the following diagram commutes
$$\xymatrix{
A_1 \ar[r]^{f_1}\ar[d]^{\varphi_1} & A_2 \ar[r]^{f_2}\ar[d]^{\varphi_2} & A_3 \ar[r]^{f_3}\ar[d]^{\varphi_3} & \cdots \ar[r]^{f_{n-1}}& A_n \ar[r]^{f_n}\ar[d]^{\varphi_n} & \Sigma A_1 \ar[d]^{\Sigma \varphi_1}\\
B_1 \ar[r]^{g_1} & B_2 \ar[r]^{g_2} & B_3 \ar[r]^{g_3} & \cdots \ar[r]^{g_{n-1}} & B_n \ar[r]^{g_n}& \Sigma B_1\\
}$$
where each row is an $n$-$\Sigma$-sequence. It is an {\em isomorphism} if $\varphi_1, \varphi_2, \varphi_3, \cdots, \varphi_n$ are all isomorphisms in $\mathcal{A}$.
\medskip

We recall the notion of a right $n$-angulated category from \cite[Definition 2.1]{[L1]}, but the Condition (RN4) is slightly different from that in \cite[Definition 2.1]{[L1]}.

\begin{definition}\label{d1}
A {\em right} $n$-\emph{angulated category} is a triple $(\mathcal{A}, \Sigma, \Theta)$, where $\mathcal{A}$ is an additive category, $\Sigma$ is an endofunctor of $\mathcal{A}$, and $\Theta$ is a class of $n$-$\Sigma$-sequences (whose elements are called right $n$-angles), which satisfies the following axioms:
\begin{itemize}
\item[(RN1)]
\begin{itemize}
\item[(a)] The class $\Theta$ is closed under isomorphisms and direct sums.

\item[(b)] For each object $A\in\mathcal{A}$ the trivial sequence
$$A\xrightarrow{1_A}A\rightarrow 0\rightarrow\cdots\rightarrow 0\rightarrow \Sigma A$$
belongs to $\Theta$.

\item[(c)] For each morphism $f_1:A_1\rightarrow A_2$ in $\mathcal{C}$, there exists a right $n$-angle whose first morphism is $f_1$.
\end{itemize}
\item[(RN2)] If an $n$-$\Sigma$-sequence belongs to $\Theta$, then its left rotation belongs to $\Theta$.

\item[(RN3)] Each solid commutative diagram
$$\xymatrix{
A_1 \ar[r]^{f_1}\ar[d]^{\varphi_1} & A_2 \ar[r]^{f_2}\ar[d]^{\varphi_2} & A_3 \ar[r]^{f_3}\ar@{-->}[d]^{\varphi_3} & \cdots \ar[r]^{f_{n-1}}& A_n \ar[r]^{f_n}\ar@{-->}[d]^{\varphi_n} & \Sigma A_1 \ar[d]^{\Sigma \varphi_1}\\
B_1 \ar[r]^{g_1} & B_2 \ar[r]^{g_2} & B_3 \ar[r]^{g_3} & \cdots \ar[r]^{g_{n-1}} & B_n \ar[r]^{g_n}& \Sigma B_1\\
}$$ with rows in $\Theta$ can be completed to a morphism of  $n$-$\Sigma$-sequences.

\item[(RN4)] Given the solid part of the diagram
  \begin{center}
    \begin{tikzpicture}
      \diagram{d}{2.5em}{2.5em}{
        A_1 & A_2 & A_3 & A_4 & \cdots & A_{n - 1} & A_n & \Sigma A_1\\
        A_1 & B_2 & B_3 & B_4 & \cdots & B_{n - 1} & B_n & \Sigma
        A_1\\
        A_2 & B_2 & C_3 & C_4 & \cdots & C_{n - 1} & C_n & \Sigma
        A_2\\
      };

      \path[->,midway,font=\scriptsize]
        (d-1-1) edge node[above] {$f_1$} (d-1-2)
        ([xshift=-0.1em] d-1-1.south) edge[-] ([xshift=-0.1em] d-2-1.north)
        ([xshift=0.1em] d-1-1.south) edge[-] ([xshift=0.1em] d-2-1.north)
        (d-1-2) edge node[above] {$f_2$} (d-1-3)
                     edge node[right] {$\varphi_2$} (d-2-2)
        (d-1-3) edge node[above] {$f_3$} (d-1-4)
                     edge[densely dashed] node[right] {$\varphi_3$} (d-2-3)
        (d-1-4) edge node[above] {$f_4$} (d-1-5)
                     edge[densely dashed] node[right] {$\varphi_4$} (d-2-4)
        (d-1-5) edge node[above] {$f_{n - 2}$} (d-1-6)
        (d-1-6) edge node[above] {$f_{n - 1}$} (d-1-7)
                     edge[densely dashed] node[right] {$\varphi_{n - 1}$} (d-2-6)
        (d-1-7) edge node[above] {$f_n$} (d-1-8)
                     edge[densely dashed] node[right] {$\varphi_n$} (d-2-7)
        ([xshift=-0.1em] d-1-8.south) edge[-] ([xshift=-0.1em] d-2-8.north)
        ([xshift=0.1em] d-1-8.south) edge[-] ([xshift=0.1em] d-2-8.north)
        (d-2-1) edge node[above] {$g_1$} (d-2-2)
                     edge node[right] {$f_1$} (d-3-1)
        (d-2-2) edge node[above] {$g_2$} (d-2-3)
        ([xshift=-0.1em] d-2-2.south) edge[-] ([xshift=-0.1em] d-3-2.north)
        ([xshift=0.1em] d-2-2.south) edge[-] ([xshift=0.1em] d-3-2.north)
        (d-2-3) edge node[above] {$g_3$} (d-2-4)
                     edge[densely dashed] node[right] {$\theta_3$} (d-3-3)
        (d-2-4) edge node[above] {$g_4$} (d-2-5)
                     edge[densely dashed] node[right] {$\theta_4$} (d-3-4)
        (d-2-5) edge node[above] {$g_{n - 2}$} (d-2-6)
        (d-2-6) edge node[above] {$g_{n - 1}$} (d-2-7)
                     edge[densely dashed] node[right] {$\theta_{n - 1}$} (d-3-6)
        (d-2-7) edge node[above] {$g_n$} (d-2-8)
                     edge[densely dashed] node[right] {$\theta_n$} (d-3-7)
        (d-2-8) edge node[right] {$\Sigma f_1$} (d-3-8)
        (d-3-1) edge node[above] {$\varphi_2$} (d-3-2)
        (d-3-2) edge node[above] {$h_2$} (d-3-3)
        (d-3-3) edge node[above] {$h_3$} (d-3-4)
        (d-3-4) edge node[above] {$h_4$} (d-3-5)
        (d-3-5) edge node[above] {$h_{n - 2}$} (d-3-6)
        (d-3-6) edge node[above] {$h_{n - 1}$} (d-3-7)
        (d-3-7) edge node[above] {$h_n$} (d-3-8)
        (d-1-4) edge[densely dashed,out=-102,in=30] node[pos=0.15,left] {$\psi_4$}
          (d-3-3)
        (d-1-7) edge[densely dashed,out=-102,in=30] node[pos=0.15,left] {$\psi_n$}
          (d-3-6);
    \end{tikzpicture}
  \end{center}
  with commuting squares and with rows in $\Theta$, the dotted
  morphisms exist such that each square commutes, and the
  $n$-$\Sigma$-sequence
\begin{center}
      \begin{tikzpicture}
        % Top row
        \node (A3) at (0,1.25){$A_3$};
        \node (A4B3) at (2,1.25){$A_4 \oplus B_3$};
        \node (A5B4C3) at (5,1.25){$A_5 \oplus B_4 \oplus C_3$};
        \node (A6B5C4) at (9,1.25){$A_6 \oplus B_5 \oplus C_4$};
        \node (tdots) at (12,1.25){$\cdots$};

        % Bottom row
        \node (mdots) at (0.25,0){$ $};
        \node (AnBn-1Cn-2) at (4.75,0){$A_n \oplus B_{n - 1} \oplus
          C_{n - 2}$};
        \node (BnCn-1) at (9.5,0){$B_n \oplus C_{n - 1}$};
        \node (mend) at (12.25,0){$C_n$};
        \node (mend2) at (14,0){$\Sigma A_3$};

        % Horizontal arrows
        \begin{scope}[font=\scriptsize,->,midway,above]
          % Top row
          \draw (A3) -- node{$\left[
              \begin{smallmatrix}
                f_3\\
                \varphi_3
              \end{smallmatrix}
            \right]$} (A4B3);
          \draw (A4B3) -- node{$\left[
              \begin{smallmatrix}
                -f_4 & 0\\
                \hfill \varphi_4 & -g_3\\
                \hfill \psi_4 & \hfill \theta_3
              \end{smallmatrix}
            \right]$} (A5B4C3);
          \draw (A5B4C3) -- node{$\left[
              \begin{smallmatrix}
                -f_5 & 0 & 0\\
                -\varphi_5 & -g_4 & 0\\
                \hfill \psi_5 & \hfill \theta_4 & h_3
              \end{smallmatrix}
            \right]$} (A6B5C4);
          \draw (A6B5C4) -- node{$\left[
              \begin{smallmatrix}
                -f_6 & 0 & 0\\
                \hfill \varphi_6 & -g_5 & 0\\
                \hfill \psi_6 & \hfill \theta_5 & h_4
              \end{smallmatrix}
            \right]$} (tdots);

          % Bottom row
          \draw (mdots) -- node{$\left[
              \begin{smallmatrix}
                -f_{n - 1} & 0 & 0\\
                (-1)^{n - 1} \varphi_{n - 1} & -g_{n - 2} & 0\\
                \psi_{n-1} & \theta_{n-2} & h_{n - 3}
              \end{smallmatrix}
            \right]$} (AnBn-1Cn-2);
          \draw (AnBn-1Cn-2) -- node{$\left[
              \begin{smallmatrix}
                (-1)^n \varphi_n & -g_{n - 1} & 0\\
                \psi_{n} & \theta_{n-1} & h_{n - 2}
              \end{smallmatrix}
            \right]$} (BnCn-1);
          \draw (BnCn-1) -- node{$\left[
              \begin{smallmatrix}
                \theta_{n} & h_{n - 1}
              \end{smallmatrix}
            \right]$} (mend);
          \draw (mend) -- node{$\Sigma f_2 \circ h_n$}
            (mend2);
        \end{scope}
      \end{tikzpicture}
    \end{center}
   belongs to $\Theta$.

   \end{itemize}
\end{definition}
The notion of a \emph{left $n$-angulated category} is defined dually.
\vspace{1mm}

If $\Sigma$ is an equivalence, and the condition in axiom (RN2) is necessary and sufficient, then the right
$n$-angulated category $(\A,\Sigma, \Theta)$ is an $n$-angulated category in the sense of Geiss-Keller-Oppermann \cite[Definition 1.1]{[GKO]} and in the sense of Bergh-Thaule \cite[Theorem 4.4]{[BT]}. If $(\A,\Sigma, \Theta)$ is a right $n$-angulated category, $(\A,\Omega, \Phi)$
is a left $n$-angulated category, $\Omega$ is a quasi-inverse of $\Sigma$ and $\Theta=\Phi$, then $(\A,\Sigma, \Theta)$ is an $n$-angulated category.

\section{One-sided $n$-suspended categories}
Let $\A$ be an additive category with an additive endofunctor $\Sigma\colon \A\to \A$ and $\X\subseteq\C$
two subcategories of $\A$. A morphism $f\colon A\to B$ with $A,B\in\C$ is called $\X$-\emph{monic}, if the
induced morphism $$\Hom_{\C}(f,\X)\colon \Hom_{\C}(B,\X)\to  \Hom_{\C}(A,\X)$$ is an epimorphism. An $\X$-\emph{epic}
is defined dually.
Recall that a morphism $x_M\colon M \to X_M$ is
called a left $\X$-approximation of $M$ if $x_M$ is an $\X$-monic and $X_M\in\X$.
Dually a morphism $x_N\colon X_N\to N$ called a right $\X$-approximation of $N$ if $x_N$ is an $\X$-epic and $X_N\in\X$.

A right $n$-$\Sigma$-sequence $A_1\xrightarrow{f_1}A_2\xrightarrow{f_2}A_3\xrightarrow{f_3}\cdots\xrightarrow{f_{n-1}}A_n\xrightarrow{f_n}\Sigma A_1$ in $\A$ is called a right \emph{$n$-$\C$-sequence} if $A_3,A_4,\cdots,A_n\in\C$ and $f_{i+1}$ is a weak cokernel
of $f_i$ for $i=1,2,\cdots, n-1$. Dually we can define a left $n$-$\C$-sequence.

\begin{definition}\label{d2}
A  \emph{right $n$-suspended category} is a triple $(\A,,R_n(\C,\Sigma),\X)$, where $\A$
  is an additive category, $\Sigma$ is an endofunctor of $\A$, $\X\subseteq\C$ are
two subcategories of $\A$ and $R_n(\C,\Sigma)$ is a class of right $n$-$\C$-sequence, which satisfies the following axioms:
\begin{itemize}
\item[(RSN1)]
\begin{itemize}
\item[(i)] The class $R_n(\C,\Sigma)$  is closed under isomorphisms and direct sums.
\item[(ii)] For each object $A\in\C$, there exists a right
$n$-$\C$-sequence
$$A\xrightarrow{x}X_2\xrightarrow{x_2}X_3\xrightarrow{x_3}
\cdots\xrightarrow{x_{n-2}}X_{n-1}\xrightarrow{x_{n-1}}B\xrightarrow{y}\Sigma A$$
where $x$ is a left $\X$-approximation of $A$ and $X_2,X_3,\cdots, X_{n-1}\in\X$.

\item[(iii)] If $A\xrightarrow{x}X_2\xrightarrow{x_2}X_3\xrightarrow{x_3}
\cdots\xrightarrow{x_{n-2}}X_{n-1}\xrightarrow{x_{n-1}}A_{n}\xrightarrow{y}\Sigma A$
is a  right $n$-$\C$-sequence, where $x$ is a left $\X$-approximation of $A$ and $X_2,X_3,\cdots, X_{n-1}\in\X$, then
for each morphism $f\colon A\to B$ in $\C$, there exists a right
$n$-$\C$-sequence
$$A\xrightarrow{\left[
              \begin{smallmatrix}
                x \\f
              \end{smallmatrix}
            \right]}X_2\oplus B\xrightarrow{~}C_3\xrightarrow{~}
\cdots\xrightarrow{~}C_{n-1}\xrightarrow{~}C_n\xrightarrow{~}\Sigma A.$$

\end{itemize}
\item[(RSN2)] For each commutative diagram of right
$n$-$\C$-sequences
$$\xymatrix{
A_1 \ar[r]^{f_1}\ar[d]^{a} & A_2 \ar[r]^{f_2}\ar[d]^{a_2} & A_3 \ar[r]^{f_3}\ar[d]^{a_3} & \cdots \ar[r]^{f_{n-2}} & A_{n-1} \ar[r]^{f_{n-1}}\ar[d]^{a_{n-1}}& A_n \ar[r]^{f_n}\ar[d]^{b} & \Sigma A_1 \ar[d]^{\Sigma a}\\
B \ar[r]^{x_1} & X_2 \ar[r]^{x_2} & X_3 \ar[r]^{x_3} & \cdots \ar[r]^{x_{n-2}} & X_{n-1} \ar[r]^{x_{n-1}} & C \ar[r]^{x_n} & \Sigma B}$$
with $X_2,X_3,\cdots,X_{n-1}\in\X$, if $a$ factors through $f_1$, then $b$ factors through $x_{n-1}$.

\item[(RSN3)] Each solid commutative diagram
$$\xymatrix{
A_1 \ar[r]^{f_1}\ar[d]^{\varphi_1} & A_2 \ar[r]^{f_2}\ar[d]^{\varphi_2} & A_3 \ar[r]^{f_3}\ar@{-->}[d]^{\varphi_3} & \cdots \ar[r]^{f_{n-1}}& A_n \ar[r]^{f_n}\ar@{-->}[d]^{\varphi_n} & \Sigma A_1 \ar[d]^{\Sigma \varphi_1}\\
B_1 \ar[r]^{g_1} & B_2 \ar[r]^{g_2} & B_3 \ar[r]^{g_3} & \cdots \ar[r]^{g_{n-1}} & B_n \ar[r]^{g_n}& \Sigma B_1\\
}$$ with rows in $R_n(\C,\Sigma)$ can be completed to a morphism of right
$n$-$\C$-sequence.

\item[(RSN4)] Given the solid part of the diagram
  \begin{center}
    \begin{tikzpicture}
      \diagram{d}{2.5em}{2.5em}{
        A_1 & A_2 & A_3 & A_4 & \cdots & A_{n - 1} & A_n & \Sigma A_1\\
        A_1 & B_2 & B_3 & B_4 & \cdots & B_{n - 1} & B_n & \Sigma
        A_1\\
        A_2 & B_2 & C_3 & C_4 & \cdots & C_{n - 1} & C_n & \Sigma
        A_2\\
      };

      \path[->,midway,font=\scriptsize]
        (d-1-1) edge node[above] {$f_1$} (d-1-2)
        ([xshift=-0.1em] d-1-1.south) edge[-] ([xshift=-0.1em] d-2-1.north)
        ([xshift=0.1em] d-1-1.south) edge[-] ([xshift=0.1em] d-2-1.north)
        (d-1-2) edge node[above] {$f_2$} (d-1-3)
                     edge node[right] {$\varphi_2$} (d-2-2)
        (d-1-3) edge node[above] {$f_3$} (d-1-4)
                     edge[densely dashed] node[right] {$\varphi_3$} (d-2-3)
        (d-1-4) edge node[above] {$f_4$} (d-1-5)
                     edge[densely dashed] node[right] {$\varphi_4$} (d-2-4)
        (d-1-5) edge node[above] {$f_{n - 2}$} (d-1-6)
        (d-1-6) edge node[above] {$f_{n - 1}$} (d-1-7)
                     edge[densely dashed] node[right] {$\varphi_{n - 1}$} (d-2-6)
        (d-1-7) edge node[above] {$f_n$} (d-1-8)
                     edge[densely dashed] node[right] {$\varphi_n$} (d-2-7)
        ([xshift=-0.1em] d-1-8.south) edge[-] ([xshift=-0.1em] d-2-8.north)
        ([xshift=0.1em] d-1-8.south) edge[-] ([xshift=0.1em] d-2-8.north)
        (d-2-1) edge node[above] {$g_1$} (d-2-2)
                     edge node[right] {$f_1$} (d-3-1)
        (d-2-2) edge node[above] {$g_2$} (d-2-3)
        ([xshift=-0.1em] d-2-2.south) edge[-] ([xshift=-0.1em] d-3-2.north)
        ([xshift=0.1em] d-2-2.south) edge[-] ([xshift=0.1em] d-3-2.north)
        (d-2-3) edge node[above] {$g_3$} (d-2-4)
                     edge[densely dashed] node[right] {$\theta_3$} (d-3-3)
        (d-2-4) edge node[above] {$g_4$} (d-2-5)
                     edge[densely dashed] node[right] {$\theta_4$} (d-3-4)
        (d-2-5) edge node[above] {$g_{n - 2}$} (d-2-6)
        (d-2-6) edge node[above] {$g_{n - 1}$} (d-2-7)
                     edge[densely dashed] node[right] {$\theta_{n - 1}$} (d-3-6)
        (d-2-7) edge node[above] {$g_n$} (d-2-8)
                     edge[densely dashed] node[right] {$\theta_n$} (d-3-7)
        (d-2-8) edge node[right] {$\Sigma f_1$} (d-3-8)
        (d-3-1) edge node[above] {$\varphi_2$} (d-3-2)
        (d-3-2) edge node[above] {$h_2$} (d-3-3)
        (d-3-3) edge node[above] {$h_3$} (d-3-4)
        (d-3-4) edge node[above] {$h_4$} (d-3-5)
        (d-3-5) edge node[above] {$h_{n - 2}$} (d-3-6)
        (d-3-6) edge node[above] {$h_{n - 1}$} (d-3-7)
        (d-3-7) edge node[above] {$h_n$} (d-3-8)
        (d-1-4) edge[densely dashed,out=-102,in=30] node[pos=0.15,left] {$\psi_4$}
          (d-3-3)
        (d-1-7) edge[densely dashed,out=-102,in=30] node[pos=0.15,left] {$\psi_n$}
          (d-3-6);
    \end{tikzpicture}
  \end{center}
  with commuting squares and with rows in $R_n(\C,\Sigma)$, and $f_{1}$, $\varphi_{2}$ are $\X$-monic, the dotted
  morphisms exist such that each square commutes, and the sequence
\begin{center}
      \begin{tikzpicture}
        % Top row
        \node (A3) at (0,1.25){$A_3$};
        \node (A4B3) at (2,1.25){$A_4 \oplus B_3$};
        \node (A5B4C3) at (5,1.25){$A_5 \oplus B_4 \oplus C_3$};
        \node (A6B5C4) at (9,1.25){$A_6 \oplus B_5 \oplus C_4$};
        \node (tdots) at (12,1.25){$\cdots$};

        % Bottom row
        \node (mdots) at (0.25,0){$ $};
        \node (AnBn-1Cn-2) at (4.75,0){$A_n \oplus B_{n - 1} \oplus
          C_{n - 2}$};
        \node (BnCn-1) at (9.5,0){$B_n \oplus C_{n - 1}$};
        \node (mend) at (12.25,0){$C_n$};
        \node (mend2) at (14,0){$\Sigma A_3$};

        % Horizontal arrows
        \begin{scope}[font=\scriptsize,->,midway,above]
          % Top row
          \draw (A3) -- node{$\left[
              \begin{smallmatrix}
                f_3\\
                \varphi_3
              \end{smallmatrix}
            \right]$} (A4B3);
          \draw (A4B3) -- node{$\left[
              \begin{smallmatrix}
                -f_4 & 0\\
                \hfill \varphi_4 & -g_3\\
                \hfill \psi_4 & \hfill \theta_3
              \end{smallmatrix}
            \right]$} (A5B4C3);
          \draw (A5B4C3) -- node{$\left[
              \begin{smallmatrix}
                -f_5 & 0 & 0\\
                -\varphi_5 & -g_4 & 0\\
                \hfill \psi_5 & \hfill \theta_4 & h_3
              \end{smallmatrix}
            \right]$} (A6B5C4);
          \draw (A6B5C4) -- node{$\left[
              \begin{smallmatrix}
                -f_6 & 0 & 0\\
                \hfill \varphi_6 & -g_5 & 0\\
                \hfill \psi_6 & \hfill \theta_5 & h_4
              \end{smallmatrix}
            \right]$} (tdots);

          % Bottom row
          \draw (mdots) -- node{$\left[
              \begin{smallmatrix}
                -f_{n - 1} & 0 & 0\\
                (-1)^{n - 1} \varphi_{n - 1} & -g_{n - 2} & 0\\
                \psi_{n-1} & \theta_{n-2} & h_{n - 3}
              \end{smallmatrix}
            \right]$} (AnBn-1Cn-2);
          \draw (AnBn-1Cn-2) -- node{$\left[
              \begin{smallmatrix}
                (-1)^n \varphi_n & -g_{n - 1} & 0\\
                \psi_{n} & \theta_{n-1} & h_{n - 2}
              \end{smallmatrix}
            \right]$} (BnCn-1);
          \draw (BnCn-1) -- node{$\left[
              \begin{smallmatrix}
                \theta_{n} & h_{n - 1}
              \end{smallmatrix}
            \right]$} (mend);
          \draw (mend) -- node{$\Sigma f_2 \circ h_n$}
            (mend2);
        \end{scope}
      \end{tikzpicture}
    \end{center}
belongs to $R_n(\C,\Sigma)$ with $\left[
              \begin{smallmatrix}
                f_3\\
                \varphi_3
              \end{smallmatrix}
            \right]$ is an $\X$-monic.
\end{itemize}
Dually, let $\A$ be an additive category with an additive endofunctor $\Omega\colon \A\to \A$ and $\X\subseteq\C$ two subcategories of $\A$.
We can define the notion of a \emph{left $n$-suspended category} $(\A,L_{n}(\C,\Omega),\X)$.
\end{definition}

\begin{remark}\label{r1}
Let $(\A,R_n(\C,\Sigma),\X)$ be a right $n$-suspended category.
For each morphism $f\colon A\to B$ in $\C$, there exists a right $n$-$\C$-sequence
$$A\xrightarrow{\left[
              \begin{smallmatrix}
                1 \\f
              \end{smallmatrix}
            \right]}A\oplus B\xrightarrow{\left[
              \begin{smallmatrix}
                f & -1
              \end{smallmatrix}
            \right]}B\xrightarrow{}0\xrightarrow{}
\cdots\xrightarrow{}0\xrightarrow{0}\Sigma A.$$
\end{remark}

\begin{remark}
In Definition \ref{d2}, when $n=3$, it is just the Definition 3.1 in \cite{[Li]}.
\end{remark}

Now we give some constructions of one-sided  $n$-suspended categories by $(n-2)$-exact categories and one-sided $n$-angulated categories.
\subsection{From $(n-2)$-exact categories to one-sided $n$-suspended categories}~\par
First, we give the first construction method  of one-sided $n$-suspended categories by  $(n-2)$-exact categories.
\begin{proposition}\label{e2}
\begin{itemize}
\item[(1)] Let $(\A,\E)$ be an $(n-2)$-exact category. Then $(\A,R_n(\A,0),\A)$
is a right $n$-suspended category, where $R_n(\A,0)=\{
A_1\to A_2\to A_3\to\cdots A_{n-1}\to A_n\to 0\mid
A_1\to A_2\to A_3\to\cdots A_{n-1}\to A_n\in\E\}$.

\item[(2)] Let $(\A,\E)$ be an $(n-2)$-exact category. Then $(\A,L_n(\A,0),\A)$
is a  left $n$-suspended category, where $L_n(\A,0)=\{
0\to A_1\to A_2\to A_3\to\cdots A_{n-1}\to A_n\mid
A_1\to A_2\to A_3\to\cdots A_{n-1}\to A_n\in\E\}$.

\item[(3)] Let $(\A,\E)$ be an $(n-2)$-exact category with enough injectives $\I$. Then $(\A,R_n(\A,0),\I)$
is a right $n$-suspended category, where $R_n(\A,0)=\{
A_1\to A_2\to A_3\to\cdots A_{n-1}\to A_n\to 0\mid
A_1\to A_2\to A_3\to\cdots A_{n-1}\to A_n\in\E\}$.

\item[(4)] Let $(\A,\E)$ be an $(n-2)$-exact category with enough projectives $\P$. Then $(\A,L_n(\A,0),\P)$
is a  left $n$-suspended category, where $L_n(\A,0)=\{
0\to A_1\to A_2\to A_3\to\cdots A_{n-1}\to A_n\mid
A_1\to A_2\to A_3\to\cdots A_{n-1}\to A_n\in\E\}$.

\end{itemize}
\end{proposition}

\proof (1) By \cite[Proposition 4.6 and Proposition 4.12]{[J]}, (RSN1) (i)  is satisfied.   (RSN1) (ii)  is trivial.
Let $$A\xrightarrow{x}X_2\xrightarrow{x_2}X_3\xrightarrow{x_3}
\cdots\xrightarrow{x_{n-2}}X_{n-1}\xrightarrow{x_{n-1}}A_{n}\rightarrow0$$
be a  right $n$-$\A$-sequence in $R_n(\A,0)$,  where $x$ is a left $\X$-approximation of $A$ and $X_2,X_3,\cdots, X_{n-1}\in\X$.  Then there is a right $n$-$\A$-sequence $A\xrightarrow{-x}X_2\xrightarrow{x_2}X_3\xrightarrow{x_3}
\cdots\xrightarrow{x_{n-2}}X_{n-1}\xrightarrow{x_{n-1}}A_{n}\rightarrow0$.
Let $f:A\rightarrow B$ be a morphism.
Since $$A\xrightarrow{-x}X_2\xrightarrow{x_2}X_3\xrightarrow{x_3}
\cdots\xrightarrow{x_{n-2}}X_{n-1}\xrightarrow{x_{n-1}}A_{n}$$ is an admissible  ($n-2$)-exact sequence, by (E2), we obtain the pushout diagram of ($-x,~x_{2},~\cdots, x_{n-2}$) along $f$
$$\xymatrix{ A\ar[r]^{-x}\ar[d]^{f}&X_{2}\ar[r]^{x_{2}}\ar[d]^{f_{2}}&X_{3}\ar[r]^{x_{3}}\ar[d]^{f_{3}}&\cdots\ar[r]^{x_{n-2}}&X_{n-1}\ar[d]^{f_{n-1}}\ar[r]^{x_{n-1}}&A_{n}\ar@{=}[d]\\
        B\ar[r]^{g_{1}}&B_{2}\ar[r]^{g_{2}}&B_{3}\ar[r]^{g_{3}}&\cdots\ar[r]^{g_{n-2}}&B_{n-1}\ar[r]^{g_{n-1}}&A_{n}}$$
By \cite[Proposition 4.8 (2)]{[J]}, we know that the mapping cone of pushout diagram
$$A\xrightarrow{\left[
              \begin{smallmatrix}
                x \\f
              \end{smallmatrix}
            \right]}X_{2}\oplus B\xrightarrow{\left[
              \begin{smallmatrix}
                -x_{2} & 0\\
                f_{2} & g_{1}
              \end{smallmatrix}
            \right]}X_{3}\oplus B_{2}\xrightarrow{\left[
              \begin{smallmatrix}
                -x_{3} & 0\\
                f_{3} & g_{2}
              \end{smallmatrix}
            \right]}X_{4}\oplus B_{3}\xrightarrow{}
\cdots\xrightarrow{}X_{n-1}\oplus B_{n-2}\xrightarrow{\left[
              \begin{smallmatrix}
                f_{n-1} & g_{n-2}
              \end{smallmatrix}
            \right]}B_{n-1}$$
is an admissible  ($n-2$)-exact sequence. Thus, we have a  right $n$-$\A$-sequence in $R_n(\A,0)$
$$A\xrightarrow{\left[
              \begin{smallmatrix}
                x \\f
              \end{smallmatrix}
            \right]}X_{2}\oplus B\xrightarrow{\left[
              \begin{smallmatrix}
                -x_{2} & 0\\
                f_{2} & g_{1}
              \end{smallmatrix}
            \right]}X_{3}\oplus B_{2}\xrightarrow{\left[
              \begin{smallmatrix}
                -x_{3} & 0\\
                f_{3} & g_{2}
              \end{smallmatrix}
            \right]}X_{4}\oplus B_{3}\xrightarrow{}
\cdots\xrightarrow{}X_{n-1}\oplus B_{n-2}\xrightarrow{\left[
              \begin{smallmatrix}
                f_{n-1} & g_{n-2}
              \end{smallmatrix}
            \right]}B_{n-1}\xrightarrow{}0.$$

 (RSN3) is trivially satisfied. Let us show (RSN2).

For each commutative diagram of right
$n$-$\A$-sequences
$$\xymatrix{
A_1 \ar[r]^{f_1}\ar[d]^{a} & A_2 \ar[r]^{f_2}\ar[d]^{a_2} & A_3 \ar[r]^{f_3}\ar[d]^{a_3} & \cdots \ar[r]^{f_{n-2}} & A_{n-1} \ar[r]^{f_{n-1}}\ar[d]^{a_{n-1}}& A_n \ar[r]\ar[d]^{b} & 0 \\
B \ar[r]^{x_1} & X_2 \ar[r]^{x_2} & X_3 \ar[r]^{x_3} & \cdots \ar[r]^{x_{n-2}} & X_{n-1} \ar[r]^{x_{n-1}} & C \ar[r] & 0}$$
with $X_2,X_3,\cdots,X_{n-1}\in\A$, assume that $a$ factors through $f_1$, that is, there exists
a morphism $u_2\colon A_2\to B$ such that $a=u_2f_1$. It follows that
$(a_2-x_1u_2)f_1=a_2f_1-x_1a=0$. Since $f_2$ is a weak cokernel of $f_1$, there exists
a morphism $u_3\colon A_3\to X_2$ such that $a_2-x_1u_2=u_3f_2$. Thus we have
$(a_3-x_2u_3)f_2=a_3f_2-x_2a_2=0$. Since $f_3$ is a weak cokernel of $f_2$, there exists
a morphism $u_4\colon A_4\to X_3$ such that $a_3-x_2u_3=u_4f_3$.
Continuing this process, there exists a morphism $u_n\colon A_n\to X_{n-1}$ such that
$a_{n-1}-x_{n-2}u_{n-1}=u_nf_{n-1}$.
It follows that $bf_{n-1}=x_{n-1}a_{n-1}=x_{n-1}(x_{n-2}u_{n-1}-u_{n}f_{1})=x_{n-1}u_nf_{n-1}$.
Since $f_{n-1}$ is an epimorphism, we have $b=x_{n-1}u_n$. This shows that $b$ factors through $x_{n-1}$.

Next we show that (RSN4).
Given the solid part of the diagram
  \begin{center}
    \begin{tikzpicture}
      \diagram{d}{2.5em}{2.5em}{
        A_1 & A_2 & A_3 & A_4 & \cdots & A_{n - 1} & A_n & 0\\
        A_1 & B_2 & B_3 & B_4 & \cdots & B_{n - 1} & B_n & 0\\
        A_2 & B_2 & C_3 & C_4 & \cdots & C_{n - 1} & C_n & 0\\
      };

      \path[->,midway,font=\scriptsize]
        (d-1-1) edge node[above] {$f_1$} (d-1-2)
        ([xshift=-0.1em] d-1-1.south) edge[-] ([xshift=-0.1em] d-2-1.north)
        ([xshift=0.1em] d-1-1.south) edge[-] ([xshift=0.1em] d-2-1.north)
        (d-1-2) edge node[above] {$f_2$} (d-1-3)
                     edge node[right] {$\varphi_2$} (d-2-2)
        (d-1-3) edge node[above] {$f_3$} (d-1-4)
                     %edge[densely dashed] node[right] {$\varphi_3$} (d-2-3)
        (d-1-4) edge node[above] {$f_4$} (d-1-5)
                     %edge[densely dashed] node[right] {$\varphi_4$} (d-2-4)
        (d-1-5) edge node[above] {$f_{n - 2}$} (d-1-6)
        (d-1-6) edge node[above] {$f_{n - 1}$} (d-1-7)
                     %edge[densely dashed] node[right] {$\varphi_{n - 1}$} (d-2-6)
        (d-1-7) edge node[above] {} (d-1-8)
                     %edge[densely dashed] node[right] {$\varphi_n$} (d-2-7)
        ([xshift=-0.1em] d-1-8.south) edge[-] ([xshift=-0.1em] d-2-8.north)
        ([xshift=0.1em] d-1-8.south) edge[-] ([xshift=0.1em] d-2-8.north)
        (d-2-1) edge node[above] {$g_1$} (d-2-2)
                     edge node[right] {$f_1$} (d-3-1)
        (d-2-2) edge node[above] {$g_2$} (d-2-3)
        ([xshift=-0.1em] d-2-2.south) edge[-] ([xshift=-0.1em] d-3-2.north)
        ([xshift=0.1em] d-2-2.south) edge[-] ([xshift=0.1em] d-3-2.north)
        (d-2-3) edge node[above] {$g_3$} (d-2-4)
                     %edge[densely dashed] node[right] {$\theta_3$} (d-3-3)
        (d-2-4) edge node[above] {$g_4$} (d-2-5)
                     %edge[densely dashed] node[right] {$\theta_4$} (d-3-4)
        (d-2-5) edge node[above] {$g_{n - 2}$} (d-2-6)
        (d-2-6) edge node[above] {$g_{n - 1}$} (d-2-7)
                     %edge[densely dashed] node[right] {$\theta_{n - 1}$} (d-3-6)
        (d-2-7) edge node[above] {} (d-2-8)
                     %edge[densely dashed] node[right] {$\theta_n$} (d-3-7)
        (d-2-8) edge node[right] {} (d-3-8)
        (d-3-1) edge node[above] {$\varphi_2$} (d-3-2)
        (d-3-2) edge node[above] {$h_2$} (d-3-3)
        (d-3-3) edge node[above] {$h_3$} (d-3-4)
        (d-3-4) edge node[above] {$h_4$} (d-3-5)
        (d-3-5) edge node[above] {$h_{n - 2}$} (d-3-6)
        (d-3-6) edge node[above] {$h_{n - 1}$} (d-3-7)
        (d-3-7) edge node[above] {} (d-3-8);
        %(d-1-4) edge[densely dashed,out=-102,in=30] node[pos=0.15,left] {$\psi_4$}
         % (d-3-3)
       % (d-1-7) edge[densely dashed,out=-102,in=30] node[pos=0.15,left] {$\psi_n$}
        %  (d-3-6)
    \end{tikzpicture}
  \end{center}
  with commuting squares and with rows in $R_n(\A,0)$, and $f_1,\varphi_2$ are $\A$-monic.
By the factorization property of weak cokernels there exist morphisms $\varphi_i\colon A_i\rightarrow B_i\ (i=3,4,\cdots, n)$, such that we have the following commutative diagram of $(n-2)$-exact sequences
$$\xymatrix{
A_1\ar[r]^{f_1}\ar@{=}[d] & A_2\ar[r]^{f_2}\ar[d]^{\varphi_2} & A_3\ar[r]^{f_3}\ar@{-->}[d]^{\varphi_3} & \cdots \ar[r]^{f_{n-2}} &  A_{n-1}\ar[r]^{f_{n-1}}\ar@{-->}[d]^{\varphi_{n-1}}  & A_{n} \ar@{-->}[d]^{\varphi_{n}}\\
A_1\ar[r]^{g_1} & B_2\ar[r]^{g_2} & B_3\ar[r]^{g_3}  & \cdots \ar[r]^{g_{n-2}} & B_{n-1}\ar[r]^{g_{n-1}} &  B_{n}.\\
}$$
By \cite[Proposition 4.8]{[J]}, we have that
$$ A_2\xrightarrow{\left[
                    \begin{smallmatrix}
                      -f_2 \\
                      \varphi_2 \\
                    \end{smallmatrix}
                  \right]}
  A_3\oplus B_2\xrightarrow{\left[
                    \begin{smallmatrix}
                      -f_3 & 0 \\
                      \varphi_3 & g_2 \\
                    \end{smallmatrix}
                  \right]}
  A_4\oplus B_3\xrightarrow{\left[
                    \begin{smallmatrix}
                      -f_4 & 0 \\
                      \varphi_4 & g_3 \\
                    \end{smallmatrix}
                  \right]}
  \cdots\xrightarrow{\left[
                    \begin{smallmatrix}
                      -f_{n-1} & 0 \\
                      \varphi_{n-1} & g_{n-2} \\
                    \end{smallmatrix}
                  \right] }
  A_{n}\oplus B_{n-1}\xrightarrow{\left[
                    \begin{smallmatrix}
                      \varphi_{n} & g_{n-1} \\
                    \end{smallmatrix}
                  \right]}
  B_{n}$$
is an $(n-2)$-exact sequence.
By the factorization property of weak cokernels there exist morphisms
 $\psi_i\colon A_i\rightarrow C_{i-1}$ $(i=4,5,\cdots, n)$ and $\theta_j\colon B_j\rightarrow C_j (j=3,4,\cdots, n)$ such that the following diagram is commutative
$$\xymatrix@R=1.5cm@C=1.2cm{
A_2\ar[r]^{\left[
                    \begin{smallmatrix}
                      -f_2 \\
                      \varphi_2 \\
                    \end{smallmatrix}
                  \right] \quad
}\ar@{=}[d] & A_3\oplus B_2\ar[r]^{\left[
                    \begin{smallmatrix}
                      -f_3 & 0 \\
                      \varphi_3 & g_2 \\
                    \end{smallmatrix}
                  \right]} \ar@{=}[d] & A_4\oplus B_3 \ar[rr]^{\ \left[
                    \begin{smallmatrix}
                      -f_4 & 0 \\
                      \varphi_4 & g_3 \\
                    \end{smallmatrix}
                  \right]}\ar[d]_{\simeq}^{\left[
                    \begin{smallmatrix}
                      -1 & 0 \\
                      0 & 1 \\
                    \end{smallmatrix}
                  \right]}& & \cdots \ar[rr]^{\left[
                    \begin{smallmatrix}
                      -f_{n-1} & 0 \\
                      \varphi_{n-1} & g_{n-2} \\
                    \end{smallmatrix}
                  \right]\qquad } & & A_{n}\oplus B_{n-1}\ar[r]^{\ \ \ \left[
                    \begin{smallmatrix}
                      \varphi_{n} & g_{n-1} \\
                    \end{smallmatrix}
                  \right]}\ar[d]_{\simeq}^{\left[
                    \begin{smallmatrix}
                      (-1)^{n+1} & 0 \\
                      0 & 1 \\
                    \end{smallmatrix}
                  \right]} \ \ \ & B_{n} \ar@{=}[d]\\
A_2\ar[r]^{\left[
                    \begin{smallmatrix}
                      -f_2 \\
                      \varphi_2 \\
                    \end{smallmatrix}
                  \right]\quad
}\ar@{=}[d] & A_3\oplus B_2\ar[r]^{\left[
                    \begin{smallmatrix}
                      f_3 & 0 \\
                      \varphi_3 & g_2 \\
                    \end{smallmatrix}
                  \right]} \ar[d]^{[0\ 1]} & A_4\oplus B_3 \ar[rr]^{\ \ \ \left[
                    \begin{smallmatrix}
                      f_4 & 0 \\
                      -\varphi_4 & g_3 \\
                    \end{smallmatrix}
                  \right]}\ar@{-->}[d]^{\left[
                    \begin{smallmatrix}
                      \psi_4 & \theta_3 \\
                    \end{smallmatrix}
                  \right]}& & \cdots\ar[rr]^{\left[
                    \begin{smallmatrix}
                      f_{n-1} & 0 \\
                      (-1)^n\varphi_{n-1} & g_{n-2} \\
                    \end{smallmatrix}
                  \right]\ \ \ \ } & & A_{n}\oplus B_{n-1}\ar[r]^{\ \ \ \tiny \left[
                    \begin{smallmatrix}
                      (-1)^{n+1}\varphi_{n} & g_{n-1} \\
                    \end{smallmatrix}
                  \right]}\ar@{-->}[d]^{ \left[
                    \begin{smallmatrix}
                      \psi_{n} & \theta_{n-1} \\
                    \end{smallmatrix}
                  \right]} \ \ \ & B_{n} \ar@{-->}[d]^{\theta_{n}}\\
A_2\ar[r]^{\varphi_2} & B_2 \ar[r]^{h_2} & C_3 \ar[rr]^{h_3}& & \cdots \ar[rr]^{h_{n-2}}& & C_{n-1} \ar[r]^{h_{n-1}} & C_{n}
}$$
where the second row is an $(n-2)$-exact sequence since it is isomorphic to the first row.
By \cite[Proposition 4.8]{[J]}, we have that
$$A_3\oplus B_2 \xrightarrow{\left[
                                                                         \begin{smallmatrix}
                                                                              -f_3 & 0 \\
                                                                              -\varphi_3 & -g_2 \\
                                                                              0 & 1\\
                                                                            \end{smallmatrix}
                                                                          \right]  }
  A_4\oplus B_3\oplus B_2\xrightarrow{\left[
                                                                            \begin{smallmatrix}
                                                                              -f_4 & 0 & 0\\
                                                                              \varphi_4 & -g_3 & 0 \\
                                                                              \psi_4 & \theta_3 & h_2\\
                                                                            \end{smallmatrix}
                                                                          \right]}
A_5\oplus B_4\oplus C_3\xrightarrow{\left[
                                                                            \begin{smallmatrix}
                                                                              -f_5 & 0  & 0\\
                                                                              -\varphi_5 & -g_4 & 0 \\
                                                                              \psi_5 & \theta_4 & h_3\\
                                                                            \end{smallmatrix}
                                                                          \right]}  \cdots$$
\begin{equation}\label{t1}
\begin{array}{l}
\xrightarrow{\left[
                                                                            \begin{smallmatrix}
                                                                              -f_{n-1} & 0 & 0 \\
                                                                              (-1)^{n+1}\varphi_{n-1} & -g_{n-2} & 0 \\
                                                                              \psi_{n-1} & \theta_{n-2} & h_{n-3}\\
                                                                            \end{smallmatrix}
                                                                          \right]}
 A_{n}\oplus B_{n-1}\oplus C_{n-2}
 \xrightarrow{\left[
                                                                            \begin{smallmatrix}
                                                                            (-1)^{n+2}\varphi_{n} & -g_{n-1} & 0 \\
                                                                              \psi_{n} & \theta_{n-1} & h_{n-2}\\
                                                                            \end{smallmatrix}
                                                                          \right]} B_{n}\oplus C_{n-1}\xrightarrow{[\theta_{n}\ h_n]} C_{n}\end{array}
\end{equation}
is an $(n-2)$-exact sequence.

The following commutative diagram
$$
\xymatrix@C=2cm@R=1.2cm{ A_3 \ar[r]^{\left[
                                                                         \begin{smallmatrix}
                                                                              f_3 \\
                                                                              \varphi_3 \\
                                                                            \end{smallmatrix}
                                                                          \right] \ \ \ }
                    \ar[d]^{\left[
                    \begin{smallmatrix}
                      -1 \\
                      0  \\
                    \end{smallmatrix}
                    \right]}
&  A_4\oplus B_3\ar[r]^{\left[
                                                                            \begin{smallmatrix}
                                                                              -f_4 & 0 \\
                                                                              \varphi_4 & -g_3 \\
                                                                              \psi_4 & \theta_3\\
                                                                            \end{smallmatrix}
                                                                          \right] \ \ \ } \ar[d]^{\left[
                                                                                                    \begin{smallmatrix}
                                                                                                      1 & 0 \\
                                                                                                      0 & 1 \\
                                                                                                      0 & 0 \\
                                                                                                    \end{smallmatrix}
                                                                                                  \right]}
&  A_5\oplus B_4\oplus C_3 \ar@{=}[d]\ar[r]&\cdots  \\
A_3\oplus B_2 \ar[r]^{\left[
                        \begin{smallmatrix}
                          -f_3 & 0 \\
                          -\varphi_3 & -g_2 \\
                          0 & 1 \\
                        \end{smallmatrix}
                      \right]} \ar[d]^{\left[
                        \begin{smallmatrix}
                          -1 & 0 \\
                        \end{smallmatrix}
                      \right]}
 & A_4\oplus B_3\oplus B_2 \ar[d]^{\left[
                        \begin{smallmatrix}
                          1 & 0 & 0\\
                          0 & 1 & g_2 \\
                        \end{smallmatrix}
                      \right]}\ar[r]^{\left[
                        \begin{smallmatrix}
                          -f_4 & 0 & 0\\
                          \varphi_4 & -g_3 & 0 \\
                          \psi_4 & \theta_3 & h_2 \\
                        \end{smallmatrix}
                      \right]} & A_5\oplus B_4\oplus C_3 \ar@{=}[d] \ar[r]&\cdots \\
 A_3 \ar[r]^{\left[
                                                                         \begin{smallmatrix}
                                                                              f_3 \\
                                                                              \varphi_3 \\
                                                                            \end{smallmatrix}
                                                                          \right] \ \ \ }
&  A_4\oplus B_3\ar[r]^{\left[
                                                                            \begin{smallmatrix}
                                                                              -f_4 & 0 \\
                                                                              \varphi_4 & -g_3 \\
                                                                              \psi_4 & \theta_3\\
                                                                            \end{smallmatrix}
                                                                          \right] \ \ \ }
&  A_5\oplus B_4\oplus C_3 \ar[r]&\cdots \\
}$$
shows that $$A_3\xrightarrow{\left[
                                                                            \begin{smallmatrix}
                                                                              f_3 \\
                                                                              \varphi_3 \\
                                                                            \end{smallmatrix}
                                                                          \right]} A_4\oplus B_3\xrightarrow{\left[
                                                                            \begin{smallmatrix}
                                                                              -f_4 & 0 \\
                                                                              \varphi_4 & -g_3 \\
                                                                              \psi_4 & \theta_3\\
                                                                            \end{smallmatrix}
                                                                          \right]} A_5\oplus B_4\oplus C_3\xrightarrow{\left[
                                                                            \begin{smallmatrix}
                                                                              -f_5 & 0 & 0 \\
                                                                              -\varphi_5 & -g_4 & 0 \\
                                                                              \psi_5 & \theta_4 & h_4\\
                                                                            \end{smallmatrix}
                                                                          \right]} \cdots\xrightarrow{\left[
                                                                            \begin{smallmatrix}
                                                                              -f_{n-1} & 0 & 0 \\
                                                                              (-1)^{n+1}\varphi_{n-1} & -g_{n-2} & 0 \\
                                                                              \psi_{n-1} & \theta_{n-2} & h_{n-3}\\
                                                                            \end{smallmatrix}
                                                                          \right]}
 A_{n}\oplus B_{n-1}\oplus C_{n-2}$$
 $$\xrightarrow{\left[
                                                                            \begin{smallmatrix}
                                                                            (-1)^{n+2}\varphi_{n} & -g_{n-1} & 0 \\
                                                                              \psi_{n} & \theta_{n-1} & h_{n-2}\\
                                                                            \end{smallmatrix}
                                                                          \right]} B_{n}\oplus C_{n-1}\xrightarrow{[\theta_{n} \ h_{n-1}]} C_{n}$$
is a direct summand of (\ref{t1}), thus it is an $(n-2)$-exact sequence.

Finally,  we claim  that $\left[
              \begin{smallmatrix}
                f_3\\
                \varphi_3
              \end{smallmatrix}
            \right]$ is an $\A$-monic. Indeed, for any morphism $u\colon A_3\to X$ with $X\in\A$, since $\varphi_2$ is an $\A$-monic, there
exists a morphism $v\colon B_2\to X$ such that $v\varphi_2=uf_2$.
It follows that $vg_1=v\varphi_2f_1=uf_2f_1$=0. So there exists a morphism
$w\colon B_3\to X$ such that $v=wg_2$.
Note that $(u-w\varphi_3)f_2=v\varphi_2-wg_2\varphi_2=0$, there exists a morphism
$w'\colon A_4\to X$ such that $u-w\varphi_2=w'f_3$ and then
$u=\left[
              \begin{smallmatrix}
                w'&w
              \end{smallmatrix}
            \right]\left[
              \begin{smallmatrix}
                f_3\\
                \varphi_3
              \end{smallmatrix}
            \right]$. This shows that $\left[
              \begin{smallmatrix}
                f_3\\
                \varphi_3
              \end{smallmatrix}
            \right]$ is an $\A$-monic.

(2) It is dual of (1).

(3) It is similar to (1).

(4) It is similar to (2).\qed
\medskip

Next, we provide more general constructions method by some special subcategories of  $(n-2)$-exact categories.

Let $\C$ be an additive category and $R_n(\C,0)$ be the class of right $(n-2)$-exact sequences
$$A_1\xrightarrow{f_1}A_2\xrightarrow{f_2}A_3\xrightarrow{f_3}\cdots\xrightarrow{f_{n-1}}A_{n}\xrightarrow{~}0$$
in $\C$.

\begin{lemma} \label{lem1}
Let $$\xymatrix{
A_1\ar[r]^{f_1}\ar@{=}[d] & A_2\ar[r]^{f_2}\ar[d]^{h_2} & A_3\ar[r]^{f_3}\ar[d]^{h_3} & \cdots \ar[r]^{f_{n-2}} & A_{n-1}\ar[r]^{f_{n-1}}\ar[d]^{h_{n-1}}  & A_{n} \ar[d]^{h_{n}}\ar[r]&0\\
A_1\ar[r]^{g_1} & B_2\ar[r]^{g_2} & B_3\ar[r]^{g_3}  & \cdots \ar[r]^{g_{n-2}} & B_{n-1}\ar[r]^{g_{n-1}} &  B_{n}\ar[r]&0
}$$ be a commutative diagram of right $(n-2)$-exact sequences in $\C$. Then
$$\xymatrix{
 A_2\ar[r]^{f_2}\ar[d]^{h_2} & A_3\ar[r]^{f_3}\ar[d]^{h_3} & \cdots \ar[r]^{f_{n-2}} & A_{n-1}\ar[r]^{f_{n-1}}\ar[d]^{h_{n-1}}  & A_{n} \ar[d]^{h_{n}}\\
 B_2\ar[r]^{g_1} & B_3\ar[r]^{g_2}  & \cdots \ar[r]^{g_{n-2}} & B_{n-1}\ar[r]^{g_{n-1}} &  B_{n}\\
}$$
is an $(n-2)$-pushout diagram, that is
$$A_2\xrightarrow{\left[
                                            \begin{smallmatrix}
                                              -f_2 \\
                                              h_2 \\
                                            \end{smallmatrix}
                                          \right]} A_3\oplus B_2\xrightarrow{\left[
                                            \begin{smallmatrix}
                                              -f_{3} & 0 \\
                                              h_{3} & g_2\\
                                            \end{smallmatrix}
                                          \right]} A_4\oplus B_3\xrightarrow{\left[
                                            \begin{smallmatrix}
                                              -f_{4} & 0 \\
                                              h_{4} & g_3 \\
                                            \end{smallmatrix}
                                          \right]}\cdots\xrightarrow{\left[
                                            \begin{smallmatrix}
                                              -f_{n-1} & 0 \\
                                              h_{n-1} & g_{n-2} \\
                                            \end{smallmatrix}
                                          \right]}A_{n}\oplus B_{n-1}\xrightarrow{\left[
                                            \begin{smallmatrix}
                                            h_{n} & g_{n-1} \\
                                            \end{smallmatrix}
                                          \right]}B_{n}\xrightarrow{~}0$$
is a right $(n-2)$-exact sequence in $\C$.
\end{lemma}

\proof See the proof of Lemma 2.8 in \cite{[L1]}. \qed

\begin{proposition}\label{p1}
Let $\X$ be an additive subcategory of\ \ $\C$. Assume the following conditions hold:
\begin{itemize}
\item[(a)] For each object $A\in\C$, there exists a sequence
$$A\xrightarrow{x}X_2\xrightarrow{x_2}X_3\xrightarrow{x_3}
\cdots\xrightarrow{x_{n-2}}X_{n-1}\xrightarrow{x_{n-1}}C\xrightarrow{~}0$$
in $R_n(\C,0)$, where $x$ is a left $\X$-approximation of $A$ and $X_2,X_3,\cdots, X_{n-1}\in\X$.
\item[(b)] If $A\xrightarrow{x}X_2\xrightarrow{x_2}X_3\xrightarrow{x_3}
\cdots\xrightarrow{x_{n-2}}X_{n-1}\xrightarrow{x_{n-1}}C\xrightarrow{~}0$
is in $R_n(\C,0)$, where $x$ is a left $\X$-approximation of $A$ and $X_2,X_3,\cdots, X_{n-1}\in\X$, then
for each morphism $f\colon A\to B$ in $\C$, there exists a sequence
$$A\xrightarrow{\left[
              \begin{smallmatrix}
                x \\f
              \end{smallmatrix}
            \right]}X_2\oplus B\xrightarrow{~}C_3\xrightarrow{~}
\cdots\xrightarrow{~}C_{n-1}\xrightarrow{~}C_n\xrightarrow{~}0$$
in $R_n(\C,0)$.
\end{itemize}
Then $(\C,R_n(\C,0),\X)$ is a  right $n$-suspended category.
\end{proposition}

\proof (RSN1) is trivially satisfied.
(RSN2) is similar to the proof of Proposition \ref{e2}(1).
(RSN3) follows from the universal property of $n$-cokernels.
(RSN4) By Lemma \ref{lem1} and combining the proof of Proposition \ref{e2}(1). \qed

\medskip

Dually, let $\C$ is an additive category and $L_n(\C,0)$ be the class of left $(n-2)$-exact sequences
$$0\xrightarrow{~}A_1\xrightarrow{f_1}A_2\xrightarrow{f_2}A_3\xrightarrow{f_3}\cdots\xrightarrow{f_{n-1}}A_{n}$$
in $\C$. We have the following.
\begin{proposition}\label{p2}
Let $\X$ be an additive subcategory of\ \ $\C$. Assume the following conditions hold:
\begin{itemize}
\item[(c)] For each object $A\in\C$, there exists a sequence
$$0\xrightarrow{~}C\xrightarrow{y_1}Y_2\xrightarrow{y_2}Y_3\xrightarrow{y_3}
\cdots\xrightarrow{y_{n-2}}Y_{n-1}\xrightarrow{y}A$$
in $L_n(\C,0)$, where $y$ is a right $\X$-approximation of $A$ and $Y_2,Y_3,\cdots, Y_{n-1}\in\X$.
\item[(d)] If $0\xrightarrow{~}C\xrightarrow{y_1}Y_2\xrightarrow{y_2}Y_3\xrightarrow{y_3}
\cdots\xrightarrow{y_{n-2}}Y_{n-1}\xrightarrow{y}A$
is in $L_n(\C,0)$, where $y$ is a right $\X$-approximation of $A$ and $Y_2,Y_3,\cdots, Y_{n-1}\in\X$, then
for each morphism $f\colon B\to A$ in $\C$, there exists a sequence
$$0\xrightarrow{~}C_1\xrightarrow{~}C_2\xrightarrow{~}C_3\xrightarrow{~}
\cdots\xrightarrow{~}C_{n-2}\xrightarrow{~}B\oplus Y_{n-1}\xrightarrow{\left[
              \begin{smallmatrix}
                f&y
              \end{smallmatrix}
            \right]}A$$
in $R_n(\C)$.
\end{itemize}
Then $(\C,L_n(\C,0),\X)$ is a left $n$-suspended category.
\end{proposition}

\begin{definition}
Let $(\A, \E)$ be an $(n-2)$-exact category and $\X\subseteq\C$ two additive subcategories of $\A$.
The subcategory $\C$ is called a \emph{special $\X$-monic closed} if it satisfies the following
conditions:
\begin{itemize}
\item[(a)] For each object $A\in\C$, there exists an admissible $(n-2)$-exact sequence
$$A\xrightarrow{x}X_2\xrightarrow{x_2}X_3\xrightarrow{x_3}
\cdots\xrightarrow{x_{n-2}}X_{n-1}\xrightarrow{x_{n-1}}C,$$
where $x$ is a left $\X$-approximation of $A$, $X_2,X_3,\cdots, X_{n-1}\in\X$ and $C\in\C$.

\item[(b)] If $A\xrightarrow{x}X_2\xrightarrow{x_2}X_3\xrightarrow{x_3}
\cdots\xrightarrow{x_{n-2}}X_{n-1}\xrightarrow{x_{n-1}}C$
is an admissible $(n-2)$-exact sequence in $\C$, where $x$ is a left $\X$-approximation of $A$ and $X_2,X_3,\cdots, X_{n-1}\in\X$, then
for each morphism $f\colon A\to B$ in $\C$, there exists an admissible $(n-2)$-exact sequence
$$A\xrightarrow{\left[
              \begin{smallmatrix}
                x \\f
              \end{smallmatrix}
            \right]}X_2\oplus B\xrightarrow{~}C_3\xrightarrow{~}
\cdots\xrightarrow{~}C_{n-1}\xrightarrow{~}C_n$$
where $C_3,C_4,\cdots,C_n\in\C$.

\end{itemize}
Dually, we can define the notion of a \emph{special $\X$-epic closed} additive subcategory of $\A$.
\end{definition}

\begin{proposition}\label{p3}
Let $(\A, \E)$ be an $(n-2)$-exact category and $\X\subseteq\C$ two additive subcategories of $\A$.
\begin{itemize}
\item[(1)] If $\C$ is a special $\X$-monic closed, then $(\A,R_n(\C,0),\X)$ is a  right $n$-angulated category, where $R_{n}(\C,0)=\{A\xrightarrow{f}B\to C_3\to \cdots \to C_n \to 0 \mid A\xrightarrow{f}B\to C_3\to \cdots \to C_n\in\E \;\textrm{and}\; C_3, \cdots,C_n\in\C\}$.

    \item[(2)]  If $\C$ is a special $\X$-epic closed, then  $(\A,L_n(\C,0),\X)$ is a left $n$-angulated category, where $L_{n}(\C,0)=\{0\to C_1\to \cdots \to C_{n-2}\to B\xrightarrow{f}A\mid C_1\to \cdots \to C_{n-2}\to B\xrightarrow{f}A\in\E \;\textrm{and}\; C_1, \cdots,C_{n-2}\in\C\}$.
\end{itemize}
\end{proposition}

\proof Using similar arguments as in the proof of Proposition \ref{p1}, we can show that the statement (1) holds by using \cite[Proposition 4.8]{[J]}. The statement (2) can be showed dually.  \qed

\begin{remark} If one takes $\C=\A$ and $\X=\A$, then Proposition \ref{p3} is just  Proposition \ref{e2} (1) and (2). If one takes $\C=\A$ with enough injectives or projectives and $\X=\P$ or $\I$, then Proposition \ref{p3} is just  Proposition \ref{e2} (3) and (4).
\end{remark}

\subsection{From one-sided  $n$-angulated categories to one-sided $n$-suspended categories}~\par
In this part, we construct right $n$-suspended categories via right $n$-angulated categories.

%\proof (1) (RSN1), (RSN2) and (RSN3) are trivially satisfied.
%Let us show (RSN4). Keeping the notations in (RSN4), by (RN4), we only need to show that $\left[
%              \begin{smallmatrix}
%                f_3\\
%                \varphi_3
%              \end{smallmatrix}
%            \right]$ is an $\X$-monic. \textcolor{red}{In fact, since $\X=0$, $\left[
%              \begin{smallmatrix}
%                f_3\\
%                \varphi_3
%              \end{smallmatrix}
%            \right]$ is an $\X$-monic.}
%(2) It is dual of (1). \qed

\begin{definition}\label{d11}
Let $(\A,\Sigma,\Theta)$ be a right $n$-angulated category and $\X\subseteq\C$ two additive subcategories of $\A$.
The subcategory $\C$ is called a \emph{special $\X$-monic closed} if it satisfies the following
conditions:
\begin{itemize}
\item[(a)] For each object $A\in\C$, there exists a right $n$-angle
$$A\xrightarrow{x}X_2\xrightarrow{x_2}X_3\xrightarrow{x_3}
\cdots\xrightarrow{x_{n-2}}X_{n-1}\xrightarrow{x_{n-1}}C\xrightarrow{x_{n}}\Sigma A,$$
where $x$ is a left $\X$-approximation of $A$,  $X_2,X_3,\cdots, X_{n-1}\in\X$ and $C\in\C$.

\item[(b)] If $A\xrightarrow{x}X_2\xrightarrow{x_2}X_3\xrightarrow{x_3}
\cdots\xrightarrow{x_{n-2}}X_{n-1}\xrightarrow{x_{n-1}}C\xrightarrow{x_{n}}\Sigma A$
is a right $n$-angle in $\C$, where $x$ is a left $\X$-approximation of $A$ and $X_2,X_3,\cdots, X_{n-1}\in\X$, then
for each morphism $f\colon A\to B$ in $\C$, there exists a right $n$-angle
$$A\xrightarrow{\left[
              \begin{smallmatrix}
                x \\f
              \end{smallmatrix}
            \right]}X_2\oplus B\xrightarrow{~}C_3\xrightarrow{~}
\cdots\xrightarrow{~}C_{n-1}\xrightarrow{~}C_n\xrightarrow{~}\Sigma A$$
where $C_3,C_4,\cdots,C_n\in\C$.
\item[(c)] If $A\xrightarrow{x}X_2\xrightarrow{x_2}X_3\xrightarrow{x_3}
\cdots\xrightarrow{x_{n-2}}X_{n-1}\xrightarrow{x_{n-1}}C\xrightarrow{x_{n}}\Sigma A$
is a right $n$-angle, where $X_2,X_3,\cdots, X_{n-1}\in\X$ and $C\in\C$, then
$x_{n-1}$ is a weak kernel of $x_{n}$.
\end{itemize}
Dually, we can define the notion of a \emph{special $\X$-epic closed} additive subcategory of $\A$.
\end{definition}
\begin{remark}\label{re3-1}
Let  $(\A,\Sigma,\Theta)$ be a right $n$-angulated category. Clearly, $\A$ is special $\X$-monic closed, whenever $\X=0$.
\end{remark}
\begin{proposition}\label{p4}
Let $(\A,\Sigma,\Theta)$ be a right $n$-angulated category and $\X\subseteq\C$ two additive subcategories of $\A$. If $\C$ is a special $\X$-monic closed, then $(\A,R_n(\C,\Sigma),\X)$ is a right $n$-suspended category, where $$R_{n}(\C,\Sigma)=\{A\xrightarrow{f}B\to C_3\to \cdots \to C_n \to \Sigma A\in\Theta \mid  C_3, \cdots,C_n\in\C\}.$$
\end{proposition}
\proof (RSN1) is trivially satisfied.
(RSN2) For each commutative diagram of right
$n$-$\C$-sequences
$$\xymatrix{
A_1 \ar[r]^{f_1}\ar[d]^{a} & A_2 \ar[r]^{f_2}\ar[d]^{a_2} & A_3 \ar[r]^{f_3}\ar[d]^{a_3} & \cdots \ar[r]^{f_{n-2}} & A_{n-1} \ar[r]^{f_{n-1}}\ar[d]^{a_{n-1}}& A_n \ar[r]^{f_n}\ar[d]^{b} & \Sigma A_1\ar[d]^{\Sigma a} \\
B \ar[r]^{x_1} & X_2 \ar[r]^{x_2} & X_3 \ar[r]^{x_3} & \cdots \ar[r]^{x_{n-2}} & X_{n-1} \ar[r]^{x_{n-1}} & C \ar[r]^{x_n} & \Sigma B}$$
with $X_2,X_3,\cdots,X_{n-1}\in\C$, assume that $a$ factors through $f_1$, that is, there exists
a morphism $u\colon A_2\to B$ such that $a=uf_1$.
It follows that $x_nb=\Sigma a\circ f_n=\Sigma u\circ(\Sigma f_1\circ f_n)=0$.
Since $x_{n-1}$ is a weak kernel of $x_{n}$, there exists a morphism $v\colon A_n\to x_{n-1}$
such that $b=x_{n-1}v$. This shows that $b$ factors through $x_{n-1}$.

(RSN3) follows from (RN3).

(RSN4) Using similar arguments as in the proof of Proposition \ref{e2}(1). \qed

Dually, we have the following.
\begin{proposition}\label{p5}
Let $(\A,\Omega,\Phi)$ be a left $n$-angulated category and $\X\subseteq\C$ two additive subcategories of $\A$.
If $\C$ is a special $\X$-epic closed, then  $(\A,L_n(\C,\Omega),\X)$ is a left $n$-suspended category, where $$L_{n}(\C,\Omega)=\{\Omega B\to C_1\to \cdots \to C_{n-2}\to B\xrightarrow{f}A\in\Phi\mid C_1, \cdots,C_{n-2}\in\C\}.$$
\end{proposition}

Consequently, by Remark \ref{re3-1}, Proposition \ref{p4} and \ref{p5}, we obtain the following result.
\begin{corollary}
\begin{itemize}
\item[(1)] If $(\A,\Sigma,\Theta)$ is a right $n$-angulated category, then
$(\A,R_n(\A,\Sigma)=\Theta,0)$ is a right $n$-suspended category.
\item[(2)] If $(\A,\Omega,\Phi)$ is a left $n$-angulated category, then
$(\A,L_n(\A,\Omega)=\Phi,0)$ is a left $n$-suspended category.
\end{itemize}
\end{corollary}

Let $(\A,\Sigma,\Theta)$ be an $n$-angulated category and $\X$ a subcategory of $\A$.
Recall that $\X$ is called \emph{extension closed} \cite[Definition 3.6]{[L2]} if for each morphism
$x_n\colon X_n\to \Sigma X_1$ with $X_1,X_n\in\X$, there exists an $n$-angle
$$X_1\xrightarrow{x_1}X_2\xrightarrow{x_2}X_3\xrightarrow{x_3}
\cdots\xrightarrow{x_{n-2}}X_{n-1}\xrightarrow{x_{n-1}}X_n\xrightarrow{x_{n}}\Sigma X_1$$
with $X_2,X_3,\cdots,X_{n-1}\in\X$.
\begin{lemma}\label{lem4}
Let $(\A,\Sigma,\Theta)$ be an $n$-angulated category and $\X\subseteq\C$ two subcategories of $\A$.
\begin{itemize}
\item[(1)] If $\C$ is extension closed and for any $A\in\C$, there exists a right $n$-angle
$$A\xrightarrow{x}X_2\xrightarrow{x_2}X_3\xrightarrow{x_3}
\cdots\xrightarrow{x_{n-2}}X_{n-1}\xrightarrow{x_{n-1}}C\xrightarrow{x_{n}}\Sigma A,$$
where $x$ is a left $\X$-approximation of $A$,  $X_2,X_3,\cdots, X_{n-1}\in\X$ and $C\in\C$, then
$\C$ is a special $\X$-monic closed.

\item[(2)] If $\C$ is extension closed and for any $A\in\C$, there exists a right $n$-angle
$$D\xrightarrow{x'}X'_2\xrightarrow{x_2}X'_3\xrightarrow{x'_3}
\cdots\xrightarrow{x'_{n-2}}X'_{n-1}\xrightarrow{x'_{n-1}}A\xrightarrow{x'_{n}}\Sigma D,$$
where $x'_n$ is a right $\X$-approximation of $A$,  $X'_2,X'_3,\cdots, X'_{n-1}\in\X$ and $D\in\C$, then
$\C$ is a special $\X$-epic closed.
\end{itemize}
\end{lemma}

\proof (1) It suffices to show that Definition \ref{d11}(b) holds.

Assume that $A\xrightarrow{x}X_2\xrightarrow{x_2}X_3\xrightarrow{x_3}
\cdots\xrightarrow{x_{n-2}}X_{n-1}\xrightarrow{x_{n-1}}C\xrightarrow{x_{n}}\Sigma A$
is a right $n$-angle in $\C$, where $x$ is a left $\X$-approximation of $A$ and $X_2,X_3,\cdots, X_{n-1}\in\X$. Then there is a right $n$-angle $A\xrightarrow{-x}X_2\xrightarrow{x_2}X_3\xrightarrow{x_3}
\cdots\xrightarrow{x_{n-2}}X_{n-1}\xrightarrow{x_{n-1}}C\xrightarrow{-x_{n}}\Sigma A$ in $\C$.

Let  $f:A\rightarrow B$ be a morphism. For the morphism $C\xrightarrow{d_n:=\Sigma f\circ x_n}\Sigma B$ with $B,C\in\C$,
since $\C$ is extension closed, then there exists an $n$-angle
$$B\xrightarrow{d_1}D_2\xrightarrow{d_2}D_3\xrightarrow{d_3}
\cdots\xrightarrow{d_{n-2}}D_{n-1}\xrightarrow{d_{n-1}}C\xrightarrow{d_{n}}\Sigma B$$
with $D_2,D_3,\cdots,D_{n-1}\in\C$.
Consider the following commutative diagram
$$\xymatrix{
\Sigma^{-1}C\ar[r]^{x_0}\ar@{=}[d]&A \ar[r]^{-x}\ar[d]^{f} & X_2 \ar[r]^{x_2} & X_3 \ar[r]^{x_3} & \cdots \ar[r]^{x_{n-2}} & X_{n-1} \ar[r]^{-x_{n-1}}& C\ar@{=}[d] \\
\Sigma^{-1}C\ar[r]^{d_0}&B \ar[r]^{d_1} & D_2 \ar[r]^{d_2} & D_3 \ar[r]^{d_3} & \cdots \ar[r]^{d_{n-2}} & D_{n-1} \ar[r]^{d_{n-1}} & C }$$
where $x_0=(-1)^n\Sigma^{-1}x_n$ and $d_0=(-1)^n\Sigma^{-1}d_n$.
By \cite[Lemma 4.1]{[GKO]}, there exists $\varphi_2,\varphi_3,\cdots,\varphi_{n-1}$ which
give the following commutative diagram of $n$-angles
$$\xymatrix{
\Sigma^{-1}C\ar[r]^{x_0}\ar@{=}[d]&A \ar[r]^{-x}\ar[d]^{f} & X_2 \ar[r]^{x_2}\ar@{-->}[d]^{\varphi_2} & X_3 \ar[r]^{x_3}\ar@{-->}[d]^{\varphi_3} & \cdots \ar[r]^{x_{n-2}} & X_{n-1} \ar[r]^{-x_{n-1}}\ar@{-->}[d]^{\varphi_{n-1}}& C\ar@{=}[d]\\
\Sigma^{-1}C\ar[r]^{d_0}&B \ar[r]^{d_1} & D_2 \ar[r]^{d_2} & D_3 \ar[r]^{d_3} & \cdots \ar[r]^{d_{n-2}} & D_{n-1} \ar[r]^{d_{n-1}} & C}$$
and the sequence
$$A\xrightarrow{\left[
              \begin{smallmatrix}
                x \\f
              \end{smallmatrix}
            \right]}X_2\oplus B\xrightarrow{~}X_3\oplus C_2\xrightarrow{~}
\cdots\xrightarrow{~}X_{n-1}\oplus D_{n-2}\xrightarrow{~}D_{n-1}\xrightarrow{~}\Sigma A$$
is an $n$-angle.
We observe that $X_3\oplus C_2,\cdots,X_{n-1}\oplus D_{n-2},D_{n-1}\in\C$.
This completes the proof.

(2) It is dual of (1).  \qed

\begin{example}\label{e100}
Let $(\A,\Sigma,\Theta)$ be an $n$-angulated category and $\X\subseteq\C$ two subcategories of $\A$.
If $(\C,\C)$ is an $\X$-mutation pair in the sense of \cite[Definition 3.1]{[L2]}, that is
for any object $A\in\C$ or $B\in\C$, there exists an $n$-angle
$$A\xrightarrow{x_1}X_2\xrightarrow{x_2}X_3\xrightarrow{x_3}
\cdots\xrightarrow{x_{n-2}}X_{n-1}\xrightarrow{x_{n-1}}B\xrightarrow{x_{n}}\Sigma A,$$
where $x_1$ is a left $\X$-approximation of $A$ and $x_{n-1}$ is a right $\X$-approximation of $B$,  $X_2,X_3,\cdots, X_{n-1}\in\X$.
Moreover, if $\C$ is extension closed, by Lemma \ref{lem4}, then we have that $\C$ is a special $\X$-monic closed and $\C$ is a special $\X$-epic closed.
Thus $(\A,R_n(\C,\Sigma),\X)$ is a right $n$-suspended category
and $(\A,L_n(\C,\Sigma^{-1}),\X)$ is a  left $n$-suspended category, where
$R_n(\C,\Sigma)$ and $L_n(\C,\Sigma^{-1})$ are constructed in Proposition \ref{p4} and Proposition \ref{p5}.
\end{example}

\section{From one-sided $n$-suspended categories to one-sided $n$-angulated
categories}
\setcounter{equation}{0}

In this section, we fix a right $n$-suspended category $(\A,R_n(\C,\Sigma),\X)$ and construct the
right $n$-angulated category on the quotient category $\C/\X$.
\medskip

For any object $A_1\in \C$, by (RSN1)(ii), there exists right
$n$-$\C$-sequences
$$\xymatrix@C=0.6cm{A_1\ar[r]^{f_1}&X_2\ar[r]^{f_2}&X_3\ar[r]^{f_3}
&\cdot\cdot\cdot\ar[r]&X_{n-1}\ar[r]^{f_{n-1}}&\mathbb{G}A_1\ar[r]^{f_n}&\Sigma A_1}$$
in $R_{n}(\C,\Sigma)$, where $f_1$ is a left $\X$-approximation of $A_1$. By Lemma \ref{y7}, $\mathbb{G}A_1$ is unique
up to isomorphism in the quotient category $\C/\X$. So for any $A_1\in\C$, we fix a right
$\C$-$n$-angles as above. For any morphism $f\in\C(A_1,B_1)$, there
exist $\varphi_i\;(i=2,3,\cdot\cdot\cdot,n-1)$ and $h$ which make the following diagram commutative.
$$\xymatrix@C=0.7cm{
A_1 \ar[r]^{f_1}\ar[d]^{f} & X_2 \ar[r]^{f_2}\ar[d]^{\varphi_2} & X_3 \ar[r]^{f_3}\ar[d]^{\varphi_3} & \cdots \ar[r]^{f_{n-2}}&X_{n-1}\ar[r]^{f_{n-1}}\ar[d]^{^{\varphi_{n-1}}}& \mathbb{G}A_1 \ar[r]^{f_n}\ar[d]^{h} & \Sigma A_1 \ar[d]^{\Sigma f}\\
B_1 \ar[r]^{g_1} & X'_2 \ar[r]^{g_2} & X'_3 \ar[r]^{g_3} & \cdots \ar[r]^{g_{n-2}}&X'_{n-1}\ar[r]^{g_{n-1}} & \mathbb{G}B_1 \ar[r]^{g_n}& \Sigma B_1\\
}$$
We define a functor $\mathbb{G}\colon\C/\X\rightarrow\C/\X$ by $\mathbb{G}(A_1)=\mathbb{G}A_1$ on the objects $A_1$ of $\C/\X$ and by $\mathbb{G}\underline{f}=\underline{h}$ on the morphisms $\underline{f}\colon A_1\rightarrow B_1$ of $\C/\X$. By Lemma \ref{y6}, $\mathbb{G}\underline{f}$ is well
defined and $\mathbb{G}$ is an additive endofunctor.

Assume that $\xymatrix@C=0.6cm{A_1\ar[r]^{f_1}&A_2\ar[r]^{f_2}&A_3\ar[r]^{f_3}
&\cdot\cdot\cdot\ar[r]^{f_{n-1}}&A_{n}\ar[r]^{f_n}&\Sigma A_1}$ is a right
$n$-$\C$-sequence in $R_n(\C,\Sigma)$, where $f_1$ is an $\X$-monic. Thus we have the following
commutative diagram
$$\xymatrix@C=0.7cm{A_1\ar[r]^{f_1}\ar@{=}[d]&A_2\ar[r]^{f_2}\ar[d]^{\varphi_2}&A_3\ar[r]^{f_3}\ar[d]^{\varphi_3}
&\cdot\cdot\cdot\ar[r]^{f_{n-2}}&A_{n-1}\ar[r]^{f_{n-1}}\ar[d]^{\varphi_{n-1}}&A_n\ar[r]^{f_n}\ar[d]^{h}&\Sigma A_1\ar@{=}[d]\\
A_1\ar[r]^{\alpha_1}&X_2\ar[r]^{\alpha_2}&X_3\ar[r]^{\alpha_3}&\cdot\cdot\cdot\ar[r]^{\alpha_{n-2}}&X_{n-1}\ar[r]^{\alpha_{n-1}}&\mathbb{G} A_1\ar[r]^{\alpha_n}&\Sigma A_1.}$$
Then we have the following $n$-$\mathbb{G}$-sequence in the quotient category $\C/\X$:
\begin{equation}\label{t2}
\begin{array}{l}\xymatrix@C=0.7cm{A_1\ar[r]^{\underline{f_1}}&A_2\ar[r]^{\underline{f_2}}&A_3\ar[r]^{\underline{f_3}}
&\cdot\cdot\cdot\ar[r]^{\underline{f_{n-1}}}&A_{n}\ar[r]^{\underline{h}\quad}&\mathbb{G}A_1.}\end{array}
\end{equation}
We define the \emph{standard right $n$-angles} in $\C/\X$ as the sextuples which are isomorphic to (\ref{t2}).
We use $\Theta'$ to denote the class of standard right $n$-angles in $\C/\X$.
\vspace{1mm}

Dually, if $(\A, L_n(\C,\Omega),\Y)$ is a left $n$-suspended category, we can construct an additive
endofunctor $\mathbb{H}\colon \C/\Y\to \C/\Y$, and define the corresponding \emph{standard left $n$-angles}.
We use $\Phi'$ to denote the class of standard left $n$-angles in $\C/\Y$.
\medskip

In order to prove our main theorem this section, we need the following lemmas.

\begin{lemma}\label{y6}
Let $$\xymatrix{
A_1 \ar[r]^{f_1}\ar[d]^{\alpha_1} & A_2 \ar[r]^{f_2}\ar[d]^{\alpha_2} & A_3 \ar[r]^{f_3}\ar[d]^{\alpha_3} & \cdots \ar[r]^{f_{n-2}} & A_{n-1} \ar[r]^{f_{n-1}}\ar[d]^{\alpha_{n-1}}& A_n \ar[r]^{f_n}\ar[d]^{\alpha_n} & \Sigma A_1 \ar[d]^{\Sigma \alpha_1}\\
A \ar[r]^{d_1} & X_2 \ar[r]^{d_2} & X_3 \ar[r]^{d_3} & \cdots \ar[r]^{d_{n-2}} & X_{n-1} \ar[r]^{d_{n-1}} & B \ar[r]^{d_n} & \Sigma A
}$$
be morphisms of right
$n$-$\C$-sequences, where $X_2,X_3,\cdots,X_{n-1}\in\X$ and $f_1$ is an $\X$-monic. If $\underline{\alpha_1}=0$, then $\underline{\alpha_n}=0$ in the quotient category $\C/\X$.
\end{lemma}

\proof By the condition $\underline{\alpha_1}=0$, there exist $X\in\X$, $x_1\colon A_1\rightarrow X$ and $x_2\colon X\rightarrow A$ such that
$\alpha_1=x_2x_1$. Since $X\in\X$ and $f_1$ is $\X$-monic, there exists $x_3\colon A_2\rightarrow X$ such that $x_3f_1=x_1$.
It follows that $\alpha_1=x_2x_1=x_2x_3f_1$, that is, $\alpha_1$ factors through $f_1$.
By (PRN2), $\alpha_n$ factors though $d_{n-1}$ implies $\underline{\alpha_n}=0$.  \qed

\begin{lemma} \label{y7}
Let $A\xrightarrow{d_1}X_2\xrightarrow{d_2}X_3\xrightarrow{d_3}\cdots\xrightarrow{d_{n-2}}X_{n-1}\xrightarrow{d_{n-1}}Y\xrightarrow{d_n}\Sigma A$ and $A\xrightarrow{d'_1}X'_2\xrightarrow{d_2}X'_3\xrightarrow{d'_3}\cdots\xrightarrow{d'_{n-2}}X'_{n-1}\xrightarrow{d'_{n-1}}Y'\xrightarrow{d'_n}\Sigma A$
be right
$\C$-$n$-angles in $R_{n}(\C,\Sigma)$, where $d_1,d'_1$ are left $\X$-approximations. Then $Y$ and $Y'$ are isomorphic in $\C/\X$.
\end{lemma}

\proof Since $d_1,d'_1$ are left $\X$-approximations, we have the following commutative
diagram.

$$\xymatrix{
A\ar[r]^{d_1}\ar@{=}[d] & X_1 \ar[r]^{d_2}\ar[d]^{\alpha_1} & X_2 \ar[r]^{d_3}\ar[d]^{\alpha_2} & \cdots \ar[r]^{d_{n-2}} & X_{n-1} \ar[r]^{d_{n-1}}\ar[d]^{\alpha_{n-1}}& Y\ar[r]^{d_n}\ar[d]^{\alpha_{n}} & \Sigma A \ar@{=}[d]\\
A \ar[r]^{d'_1}\ar@{=}[d] & X'_1 \ar[r]^{d'_2}\ar[d]^{\alpha_1'} & X'_2 \ar[r]^{d'_3}\ar[d]^{\alpha_2'} & \cdots \ar[r]^{d'_{n-2}} & X'_{n-1} \ar[r]^{d'_{n-1}}\ar[d]^{\alpha_{n-1}'}& Y' \ar[r]^{d'_n}\ar[d]^{\alpha'_{n}} & \Sigma A\ar@{=}[d]\\
A\ar[r]^{d_1} & X_1 \ar[r]^{d_2} & X_2\ar[r]^{d_3} & \cdots \ar[r]^{d_{n-2}} & X_{n-1}  \ar[r]^{d_{n-1}} & Y \ar[r]^{d_n}& \Sigma A}$$
By Lemma \ref{y6}, we have $\overline{\alpha'_n}\circ \overline{\alpha_n}=\textrm{Id}_{\overline{Y}}$.
Similarly, we can prove that $\overline{\alpha_n}\circ \overline{\alpha'_n}=\textrm{Id}_{\overline{Y'}}$.
Therefore, $Y$ and $Y'$ are isomorphic in $\C/\X$.  \qed
\begin{lemma}\label{y8}
Assume that we have a commutative diagram
$$\xymatrix{
A_1 \ar[r]^{f_1}\ar[d]^{\varphi_1} & A_2 \ar[r]^{f_2}\ar[d]^{\varphi_2} & A_3 \ar[r]^{f_3}\ar[d]^{\varphi_3} & \cdots \ar[r]^{f_{n-1}}& A_n \ar[r]^{f_n}\ar[d]^{\varphi_n} & \Sigma A_1 \ar[d]^{\Sigma \varphi_1}\\
B_1 \ar[r]^{g_1} & B_2 \ar[r]^{g_2} & B_3 \ar[r]^{g_3} & \cdots \ar[r]^{g_{n-1}} & B_n \ar[r]^{g_n}& \Sigma B_1\\
}$$
where the rows are right $n$-$\C$-sequences in $R_n(\C,\Sigma)$ and $f_1, g_1$ are $\X$-monic. Then we have the following commutative diagram
$$\xymatrix{
A_1 \ar[r]^{\underline{f_1}}\ar[d]^{\underline{\varphi_1}} & A_2 \ar[r]^{\underline{f_2}}\ar[d]^{\underline{\varphi_2}} & A_3 \ar[r]^{\underline{f_3}}\ar[d]^{\underline{\varphi_3}} & \cdots \ar[r]^{\underline{f_{n-1}}}& A_n \ar[r]^{\underline{\alpha_n}}\ar[d]^{\underline{\varphi_n}} & \mathbb{G}A_1 \ar[d]^{\mathbb{G}\underline{\varphi_1}}\\
B_1 \ar[r]^{\underline{g_1}} & B_2 \ar[r]^{\underline{g_2}} & B_3 \ar[r]^{\underline{g_3}} & \cdots \ar[r]^{\underline{g_{n-1}}} & B_n \ar[r]^{\underline{\beta_n}}& \mathbb{G}B_1
}$$ where the rows are standard right $n$-angles in $\C/\X$.
\end{lemma}

\proof Consider the following  commutative diagrams where three rows are right $n$-$\C$-sequences in $R_n(\C,\Sigma)$, and $\mathbb{G}\underline{\varphi_1}=\underline{\theta_1}$:
$$\xymatrix{
A_1 \ar[r]^{f_1}\ar@{=}[d] & A_2 \ar[r]^{f_2}\ar[d]^{\alpha_2} & A_3 \ar[r]^{f_3}\ar[d]^{\alpha_3} & \cdots \ar[r]^{f_{n-2}} & A_{n-1} \ar[r]^{f_{n-1}}\ar[d]^{\alpha_{n-1}}& A_n \ar[r]^{f_n}\ar[d]^{\alpha_n} & \Sigma A_1 \ar@{=}[d]\\
A_1 \ar[r]^{d_1}\ar[d]^{\varphi_1} & X_2 \ar[r]^{d_2}\ar[d]^{\alpha_1'} & X_3 \ar[r]^{d_3}\ar[d]^{\alpha_2'} & \cdots \ar[r]^{d_{n-2}} & X_{n-1} \ar[r]^{d_{n-1}}\ar[d]^{\alpha_{n-2}'}& \mathbb{G}A_1 \ar[r]^{d_n}\ar[d]^{\theta_1} & \Sigma A_1\ar[d]^{\Sigma \varphi_1}\\
B_1 \ar[r]^{d_1'} & X_2' \ar[r]^{d_2'} & X_3' \ar[r]^{d_3'} & \cdots \ar[r]^{d_{n-2}'} & X_{n-1}'  \ar[r]^{d_{n-1}'} & \mathbb{G}B_1 \ar[r]^{d_n'}& \Sigma B_1
}$$
and
$$\xymatrix{
A_1 \ar[r]^{f_1}\ar[d]^{\varphi_1} & A_2 \ar[r]^{f_2}\ar[d]^{\varphi_2} & A_3 \ar[r]^{f_3}\ar[d]^{\varphi_3} & \cdots\ar[r]^{f_{n-2}} & A_{n-1}  \ar[r]^{f_{n-1}}\ar[d]^{\varphi_{n-1}}& A_n \ar[r]^{f_n}\ar[d]^{\varphi_n} & \Sigma A_1 \ar[d]^{\Sigma \varphi_1}\\
B_1 \ar[r]^{g_1}\ar@{=}[d] & B_2 \ar[r]^{g_2}\ar[d]^{\beta_2} & B_3 \ar[r]^{g_3}\ar[d]^{\beta_3} & \cdots \ar[r]^{g_{n-2}} & B_{n-1} \ar[r]^{g_{n-1}}\ar[d]^{\beta_{n-1}} & B_n \ar[r]^{g_n}\ar[d]^{\beta_n}& \Sigma B_1\ar@{=}[d]\\
B_1 \ar[r]^{d_1'} & X_2' \ar[r]^{d_2'} & X_3' \ar[r]^{d_3'} & \cdots \ar[r]^{d_{n-2}'} & X_{n-1}'  \ar[r]^{d_{n-1}'} & \mathbb{G}B_1 \ar[r]^{d_n'}& \Sigma B_1
}$$
 By Lemma \ref{y6}, we have
$\mathbb{G}\underline{\varphi_1}\circ \underline{\alpha_n}=\underline{\beta_n}\circ\underline{\varphi_n}$. Then the following diagram commutes
$$\xymatrix{
A_1 \ar[r]^{\underline{f_1}}\ar[d]^{\underline{\varphi_1}} & A_2 \ar[r]^{\underline{f_2}}\ar[d]^{\underline{\varphi_2}} & A_3 \ar[r]^{\underline{f_3}}\ar[d]^{\underline{\varphi_3}} & \cdots \ar[r]^{\underline{f_{n-1}}}& A_n \ar[r]^{\underline{\alpha_n}}\ar[d]^{\underline{\varphi_n}} & \mathbb{G}A_1 \ar[d]^{\mathbb{G}\underline{\varphi_1}}\\
B_1 \ar[r]^{\underline{g_1}} & B_2 \ar[r]^{\underline{g_2}} & B_3 \ar[r]^{\underline{g_3}} & \cdots \ar[r]^{\underline{g_{n-1}}} & B_n \ar[r]^{\underline{\beta_n}}& \mathbb{G}B_1.
}$$
\qed

Now we prove our first main result in this section.

\begin{theorem}\label{y9}
\begin{itemize}
\item[(1)]If $(\A,R_n(\C,\Sigma),\X)$ is a  right $n$-suspended category, then
$(\C/\X,\mathbb{G},\Theta')$ is a right $n$-angulated category.

\item[(2)]If $(\A,L_n(\C,\Omega),\X)$ is a left $n$-suspended category, then
$(\C/\Y,\mathbb{H},\Phi')$ is a left $n$-angulated category.
\end{itemize}
\end{theorem}

\proof (1) We will check the axioms of right $n$-angulated categories.

It is easy to see that the class of standard right $n$-angles $\Theta'$ is closed under isomorphisms and direct sums. So (RN1)(a) is satisfied.

For any object $A\in\C$, the identity morphism of $A$ is an $\X$-monic. The commutative diagram
$$\xymatrix{
A \ar[r]^{1_A}\ar@{=}[d] & A \ar[r]^{0}\ar[d]^{x_1} & 0 \ar[r]^{0}\ar[d]^{0} & \cdots \ar[r]^{0} & 0 \ar[r]^{0}\ar[d]^{0}& 0 \ar[r]^{0}\ar[d]^{0} & \Sigma A \ar@{=}[d]\\
A \ar[r]^{x_1} & X_2 \ar[r]^{x_2} & X_3 \ar[r]^{x_3} & \cdots \ar[r]^{d_{n-2}} & X_{n-1} \ar[r]^{x_{n-1}}& \mathbb{G}A \ar[r]^{x_n} & \Sigma A \\
}$$ shows that $$A\xrightarrow{\underline{1_A}}A\rightarrow 0\rightarrow\cdots\rightarrow 0\rightarrow \mathbb{G}A$$ belongs to $\Theta'$. Thus (RN1)(b) is satisfied.

For any morphism $f\colon A \to B$ in $\C$, By (RSN1)(ii) and (RSN1)(iii),  there exist  right
$n$-$\C$-angles
$$A\xrightarrow{\left[
              \begin{smallmatrix}
                x \\f
              \end{smallmatrix}
            \right]}X_2\oplus B\xrightarrow{~}C_3\xrightarrow{~}
\cdots\xrightarrow{~}C_{n-1}\xrightarrow{~}C_n\xrightarrow{~}\Sigma A$$
where  $x$ is a left $\X$-approximation of $A$. Note that $\left[
              \begin{smallmatrix}
                x \\f
              \end{smallmatrix}
            \right]$
is an $\X$-monic, the above right
$n$-$\C$-sequence induces a standard right $n$-angle
$$A\xrightarrow{\underline{f}}B\xrightarrow{~}C_3\xrightarrow{~}\cdots\xrightarrow{~}
C_{n-1}\xrightarrow{~}C_n\xrightarrow{~}\mathbb{G}A.$$
Thus (RN1)(c) is satisfied.

(RN2)
Let
$A_1\xrightarrow{\underline{f_1}}A_2\xrightarrow{\underline{f_2}}A_3\xrightarrow{\underline{f_3}}\cdots\xrightarrow{\underline{f_{n-1}}}
A_n\xrightarrow{\underline{\alpha_n}} \mathbb{G}A_1$ be a standard right $n$-angles in $\C/\X$.

For each object $A_1\in\C$, By (RSN1)(ii), there exists a right
$n$-$\C$-sequence
$$A_1\xrightarrow{x}X_2\xrightarrow{x_2}X_3\xrightarrow{x_3}
\cdots\xrightarrow{x_{n-2}}X_{n-1}\xrightarrow{x_{n-1}}\mathbb{G}A_1\xrightarrow{y}\Sigma A_1$$
where $x$ is a left $\X$-approximation of $A$ and $X_2,X_3,\cdots, X_{n-1}\in\X$.
Without loss of generality, we can assume that
$$A_1\xrightarrow{\left[
              \begin{smallmatrix}
                x \\f_1
              \end{smallmatrix}
            \right]}X_2\oplus A_2\xrightarrow{\left[
              \begin{smallmatrix}
                z &f_2
              \end{smallmatrix}
            \right]}A_3\xrightarrow{f_3}
\cdots\xrightarrow{f_{n-2}}A_{n-1}\xrightarrow{f_{n-1}}A_n\xrightarrow{f_n}\Sigma A_1$$
induces the standard right $n$-angle $A_1\xrightarrow{\underline{f_1}}A_2\xrightarrow{\underline{f_2}}A_3\xrightarrow{\underline{f_3}}\cdots\xrightarrow{\underline{f_{n-1}}}
A_n\xrightarrow{\underline{\alpha_n}} \mathbb{G}A_1$, that is, there exists the following
commutative diagram
\begin{equation}\label{t3}
\begin{array}{l}\xymatrix@C=0.7cm{A_1\ar[r]^{\left[
              \begin{smallmatrix}
                x \\f_1
              \end{smallmatrix}
            \right]\quad\;\;}\ar@{=}[d]&X_2\oplus A_2\ar[r]^{\quad \left[
              \begin{smallmatrix}
                z &f_2
              \end{smallmatrix}
            \right]}\ar[d]^{\delta_2}&A_3\ar[r]^{f_3}\ar[d]^{\delta_3}
&\cdot\cdot\cdot\ar[r]^{f_{n-2}}&A_{n-1}\ar[r]^{f_{n-1}}\ar[d]^{\delta_{n-1}}&A_n\ar[r]^{f_n}\ar[d]^{\alpha_n}&\Sigma A_1\ar@{=}[d]\\
A_1\ar[r]^{x}&X_2\ar[r]^{x_2}&X_3\ar[r]^{x_3}&\cdot\cdot\cdot\ar[r]^{x_{n-2}}&X_{n-1}\ar[r]^{x_{n-1}}&\mathbb{G} A_1\ar[r]^{y}&\Sigma A_1}\end{array}
\end{equation}
of right
$n$-$\C$-sequences in $R_n(\C,\Sigma)$.

For the morphism $f_1\colon A_1\to A_2$ in $\C$, by Remark \ref{r1}, there exists an $n$-$\C$-sequence
$$A_1\xrightarrow{\left[
              \begin{smallmatrix}
                1 \\f_1
              \end{smallmatrix}
            \right]}A_1\oplus A_2\xrightarrow{\left[
              \begin{smallmatrix}
                f_1 & -1
              \end{smallmatrix}
            \right]}A_2\xrightarrow{}0\xrightarrow{}
\cdots\xrightarrow{}0\xrightarrow{0}\Sigma A_1.$$
By (RSN1)(i), we have that
$$A_1\oplus A_2\xrightarrow{\left[
              \begin{smallmatrix}
                x &0\\
                0&1
              \end{smallmatrix}
            \right]}X_2\oplus A_2\xrightarrow{\left[
              \begin{smallmatrix}
                x_2&0
              \end{smallmatrix}
            \right]}X_3\xrightarrow{x_3}
\cdots\xrightarrow{x_{n-2}}X_{n-1}\xrightarrow{x_{n-1}}\mathbb{G}A_1\xrightarrow{\left[
              \begin{smallmatrix}
                y \\0
              \end{smallmatrix}
            \right]}\Sigma A_1\oplus \Sigma A_2$$
is a right
$n$-$\C$-sequence. By (RSN4), we have the following commutative
\begin{center}
    \begin{tikzpicture}
      \diagram{d}{2.5em}{2.5em}{
        A_1 & A_1\oplus A_2 & 0 & 0 & \cdots & 0 & 0 & \Sigma A_1\\
        A_1 & X_2\oplus A_2 & A_3 & A_4 & \cdots & A_{n - 1} & A_n & \Sigma
        A_1\\
        A_1\oplus A_2 &X_2\oplus A_2 & X_3 & X_4 & \cdots & X_{n - 1} & \mathbb{G}A_1 & \Sigma
        A_1\oplus\Sigma A_2\\
      };

      \path[->,midway,font=\scriptsize]
        (d-1-1) edge node[above] {$\left[
              \begin{smallmatrix}
                1 \\f_1
              \end{smallmatrix}
            \right]$} (d-1-2)
        ([xshift=-0.1em] d-1-1.south) edge[-] ([xshift=-0.1em] d-2-1.north)
        ([xshift=0.1em] d-1-1.south) edge[-] ([xshift=0.1em] d-2-1.north)
        (d-1-2) edge node[above] {$\left[
              \begin{smallmatrix}
                f_1&-1
              \end{smallmatrix}
            \right]$} (d-1-3)
                     edge node[right] {$\left[
              \begin{smallmatrix}
                x&0\\0&1
              \end{smallmatrix}
            \right]$} (d-2-2)
        (d-1-3) edge node[above] {} (d-1-4)
                     edge[densely dashed] node[right] {$\varphi_3$} (d-2-3)
        (d-1-4) edge node[above] {} (d-1-5)
                     edge[densely dashed] node[right] {$\varphi_4$} (d-2-4)
        (d-1-5) edge node[above] {} (d-1-6)
        (d-1-6) edge node[above] {} (d-1-7)
                     edge[densely dashed] node[right] {$\varphi_{n - 1}$} (d-2-6)
        (d-1-7) edge node[above] {} (d-1-8)
                     edge[densely dashed] node[right] {$\varphi_n$} (d-2-7)
        ([xshift=-0.1em] d-1-8.south) edge[-] ([xshift=-0.1em] d-2-8.north)
        ([xshift=0.1em] d-1-8.south) edge[-] ([xshift=0.1em] d-2-8.north)
        (d-2-1) edge node[above] {$\left[
              \begin{smallmatrix}
                x\\f_1
              \end{smallmatrix}
            \right]$} (d-2-2)
                     edge node[right] {$\left[
              \begin{smallmatrix}
                1 \\f_1
              \end{smallmatrix}
            \right]$} (d-3-1)
        (d-2-2) edge node[above] {$\left[
              \begin{smallmatrix}
                z&f_2
              \end{smallmatrix}
            \right]$} (d-2-3)
        ([xshift=-0.1em] d-2-2.south) edge[-] ([xshift=-0.1em] d-3-2.north)
        ([xshift=0.1em] d-2-2.south) edge[-] ([xshift=0.1em] d-3-2.north)
        (d-2-3) edge node[above] {$f_3$} (d-2-4)
                     edge[densely dashed] node[right] {$\theta_3$} (d-3-3)
        (d-2-4) edge node[above] {$f_4$} (d-2-5)
                     edge[densely dashed] node[right] {$\theta_4$} (d-3-4)
        (d-2-5) edge node[above] {$f_{n - 2}$} (d-2-6)
        (d-2-6) edge node[above] {$f_{n - 1}$} (d-2-7)
                     edge[densely dashed] node[right] {$\theta_{n - 1}$} (d-3-6)
        (d-2-7) edge node[above] {$f_n$} (d-2-8)
                     edge[densely dashed] node[right] {$\theta_n$} (d-3-7)
        (d-2-8) edge node[right] {$\Sigma\left[
              \begin{smallmatrix}
                1 \\f_1
              \end{smallmatrix}
            \right]$} (d-3-8)
        (d-3-1) edge node[above] {$\left[
              \begin{smallmatrix}
                x&0\\0&1
              \end{smallmatrix}
            \right]$} (d-3-2)
        (d-3-2) edge node[above] {$\left[
              \begin{smallmatrix}
                x_2 &0
              \end{smallmatrix}
            \right]$} (d-3-3)
        (d-3-3) edge node[above] {$x_3$} (d-3-4)
        (d-3-4) edge node[above] {$x_4$} (d-3-5)
        (d-3-5) edge node[above] {$x_{n - 2}$} (d-3-6)
        (d-3-6) edge node[above] {$x_{n - 1}$} (d-3-7)
        (d-3-7) edge node[above] {$\left[
              \begin{smallmatrix}
                y \\0
              \end{smallmatrix}
            \right]$} (d-3-8)
        (d-1-4) edge[densely dashed,out=-102,in=30] node[pos=0.15,left] {$\psi_4$}
          (d-3-3)
        (d-1-7) edge[densely dashed,out=-102,in=30] node[pos=0.15,left] {$\psi_n$}
          (d-3-6);
    \end{tikzpicture}
  \end{center}
of right $n$-$\C$-sequences in $R_n(\C,\Sigma)$ and $$ A_2\xrightarrow{\varphi_3=-f_2}
  A_3\xrightarrow{\left[
                    \begin{smallmatrix}
                      -f_3 \\
                      \theta_3
                    \end{smallmatrix}
                  \right]}
  A_4\oplus X_3\xrightarrow{\left[
                    \begin{smallmatrix}
                      -f_4 & 0 \\
                      \theta_4 & x_3 \\
                    \end{smallmatrix}
                  \right]}
  \cdots\xrightarrow{\left[
                    \begin{smallmatrix}
                      -f_{n-1} & 0 \\
                      \theta_{n-1} & x_{n-2} \\
                    \end{smallmatrix}
                  \right] }
  A_{n}\oplus X_{n-1}\xrightarrow{\left[
                    \begin{smallmatrix}
                      \theta_{n} & x_{n-1} \\
                    \end{smallmatrix}
                  \right]}
  \mathbb{G}A_1\xrightarrow{\Sigma f_1\circ y}\Sigma A_2$$
is a right $n$-$\C$-sequence in $R_n(\C,\Sigma)$ with $\varphi_2$ is an $\X$-monic. The following commutative diagram with rows are right $n$-$\C$-sequences in $R_n(\C,\Sigma)$
$$\xymatrix@C=1.2cm@R=1.2cm{A_1\ar[r]^{x}\ar[d]^{f_1}&X_2\ar[r]^{ x_2}\ar[d]^{z}&X_3\ar[r]^{x_3}\ar[d]^{\left[
                    \begin{smallmatrix}
                      0 \\
                      1\\
                    \end{smallmatrix}
                  \right]}
&\cdots\ar[r]^{x_{n-2}}&X_{n-1}\ar[r]^{x_{n-1}}\ar[d]^{\left[
                    \begin{smallmatrix}
                      0 \\
                      1\\
                    \end{smallmatrix}
                  \right]}&\mathbb{G}A_1\ar[r]^{y}\ar@{=}[d]&\Sigma A_1\ar[d]^{\Sigma f_1}\\
A_2\ar[r]^{-f_2}\ar@{=}[d]&A_3\ar[r]^{\left[\begin{smallmatrix}
                      -f_3 \\
                      \theta_3\\
                    \end{smallmatrix}
                  \right]\quad}\ar@{-->}[d]^{\omega_2}&A_4\oplus X_3\ar[r]^{\left[\begin{smallmatrix}
                      -f_4&0 \\
                      \theta_4&x_3
                    \end{smallmatrix}
                  \right]}\ar@{-->}[d]^{\omega_3}
&\cdots\ar[r]^{\left[\begin{smallmatrix}
                      -f_{n-1}&0 \\
                      \theta_{n-1}&x_{n-2}
                    \end{smallmatrix}
                  \right]\quad}&A_{n}\oplus X_{n-1}\ar[r]^{\quad \left[
                    \begin{smallmatrix}
                      \theta_{n} & x_{n-1} \\
                    \end{smallmatrix}
                  \right]}\ar@{-->}[d]^{\omega_{n-1}}&\mathbb{G}A_1\ar[r]^{\Sigma f_1\circ y}\ar@{-->}[d]^{\omega_n}&\Sigma A_2\ar@{=}[d]\\
A_2\ar[r]^{x'}&X'_2\ar[r]^{x'_2}&X'_3\ar[r]^{x'_3}&\cdot\cdot\cdot\ar[r]^{x'_{n-2}}&X'_{n-1}\ar[r]^{x'_{n-1}}&\mathbb{G} A_2\ar[r]^{y'}&\Sigma A_2}$$
implies that $\underline{\omega_n}=\mathbb{G}\underline{f_1}$ and
$$A_2\xrightarrow{\underline{-f_2}}A_3\xrightarrow{-\underline{f_3}}A_4\xrightarrow{-\underline{f_4}}\cdots\xrightarrow{-\underline{f_{n-1}}}
A_{n}\xrightarrow{\underline{\theta_{n}}} \mathbb{G}A_1\xrightarrow{\underline{\omega_n}} \mathbb{G}A_2
$$
is a standard right $n$-angle in $\C/\X$. The following commutative diagram
$$\xymatrix@C=1cm{A_2\ar[r]^{\underline{-f_2}}\ar@{=}[d]& A_3\ar[r]^{\underline{-f_3}}\ar[d]^{-1}&A_4\ar[r]^{\underline{-f_4}}\ar@{=}[d]
&\cdot\cdot\cdot\ar[r]^{\underline{-f_{n-1}}}&A_{n}\ar[r]^{\underline{\theta_{n}}}\ar[d]^{(-1)^n}&\mathbb{G}A_1\ar[r]^{\mathbb{G}\underline{f_1}}\ar[d]^{(-1)^n}&\mathbb{G}A_2\ar@{=}[d]\\
A_2\ar[r]^{\underline{f_2}}&A_3\ar[r]^{\underline{f_3}}&A_4\ar[r]^{\underline{f_4}}&\cdot\cdot\cdot\ar[r]^{\underline{f_{n-1}}}&A_{n}\ar[r]^{\underline{\theta_{n}}}&\mathbb{G} A_1\ar[r]^{(-1)^n\mathbb{G}\underline{f_1}}&\mathbb{G} A_2}$$
shows that
$$A_2\xrightarrow{\underline{f_2}}A_3\xrightarrow{\underline{f_3}}A_4\xrightarrow{\underline{f_4}}\cdots\xrightarrow{\underline{f_{n-1}}}
A_{n}\xrightarrow{\underline{\theta_{n}}} \mathbb{G}A_1\xrightarrow{(-1)^n\mathbb{G}\underline{f_1}} \mathbb{G}A_2
$$
belongs to $\Theta'$. It remains to show that $\underline{\theta_n}=\underline{\alpha_n}$.
Consider the following commutative diagram
\begin{equation}\label{t4}
\begin{array}{l}\xymatrix@C=0.7cm{A_1\ar[r]^{\left[
              \begin{smallmatrix}
                x \\f_1
              \end{smallmatrix}
            \right]\quad\;\;}\ar@{=}[d]&X_2\oplus A_2\ar[r]^{\quad \left[
              \begin{smallmatrix}
                z &f_2
              \end{smallmatrix}
            \right]}\ar[d]^{\left[
              \begin{smallmatrix}
                1 &0
              \end{smallmatrix}
            \right]}&A_3\ar[r]^{f_3}\ar[d]^{\theta_3}
&\cdot\cdot\cdot\ar[r]^{f_{n-2}}&A_{n-1}\ar[r]^{f_{n-1}}\ar[d]^{\theta_{n-1}}&A_n\ar[r]^{f_n}\ar[d]^{\theta_n}&\Sigma A_1\ar@{=}[d]\\
A_1\ar[r]^{x}&X_2\ar[r]^{x_2}&X_3\ar[r]^{x_3}&\cdot\cdot\cdot\ar[r]^{x_{n-2}}&X_{n-1}\ar[r]^{x_{n-1}}&\mathbb{G} A_1\ar[r]^{y}&\Sigma A_1}\end{array}
\end{equation}
and the commutative diagram (\ref{t3}), by Lemma \ref{y6}, we have $\underline{\theta_n}=\underline{\alpha_n}$.
\vspace{1mm}

(RN3) Suppose that there exists a commutative diagram
\begin{equation}\label{t5}
\begin{array}{l}
\xymatrix{
A_1 \ar[r]^{\underline{f_1}}\ar[d]^{\underline{h_1}} & A_2 \ar[r]^{\underline{f_2}}\ar[d]^{\underline{h_2}} & A_3 \ar[r]^{\underline{f_3}}  & \cdots \ar[r]^{\underline{f_{n-1}}}& A_{n} \ar[r]^{\underline{\alpha_{n}}} & \mathbb{G}A_1\ar[d]^{\mathbb{G}\underline{h_1}}\\
B_1 \ar[r]^{\underline{g_1}} & B_2 \ar[r]^{\underline{g_2}} & B_3 \ar[r]^{\underline{g_3}} & \cdots \ar[r]^{\underline{g_{n-1}}} & B_{n} \ar[r]^{\underline{\beta_{n}}}&\mathbb{G}B_1\\
}\end{array}
\end{equation}
with rows standard right $(n+2)$-angles in $\C/\X$.
By assumption, the equality $\underline{h_2}\circ \underline{f_1}=\underline{g_1}\circ \underline{h_1}$ holds, there exist morphisms $u\colon A_1\to X$ and $v\colon X\to B_2$ with $X\in\X$ such that
$h_2f_1-g_1h_1=vu$. Since $f_1$ is an $\X$-monic, there exists a morphism
$w\colon A_2\to X$ such that $wf_1=u$. It follows that $(h_2-vw)f_1=h_1f_1-vu=g_1h_1$. By (RSN3), we have the
following commutative diagram
\begin{equation}\label{t6}
\begin{array}{l}
\xymatrix{
A_1 \ar[r]^{{f_1}}\ar[d]^{{h_1}} & A_2 \ar[r]^{{f_2}}\ar[d]^{{h_2-vw}} & A_3 \ar@{-->}[d]\ar[r]^{{f_3}}  & \cdots \ar[r]^{{f_{n-1}}}& A_{n}\ar@{-->}[d] \ar[r]^{{f_{n}}} & \Sigma A_1\ar[d]^{\Sigma{h_1}}\\
B_1 \ar[r]^{{g_1}} & B_2 \ar[r]^{{g_2}} & B_3 \ar[r]^{{g_3}} & \cdots \ar[r]^{{g_{n-1}}} & B_{n} \ar[r]^{{g_{n}}}&\Sigma B_1\\
}\end{array}
\end{equation}
of right $n$-$\C$-sequences in $R_n(\C,\Sigma)$. By Lemma \ref{y8} and $\underline{h_2-vw}=\underline{h_2}$, the diagram (\ref{t5}) can be completed to a
morphism of right $n$-angles follows from diagram (\ref{t6}).
\vspace{1mm}

(RN4)  Given the solid part of the diagram
\begin{center}
    \begin{tikzpicture}
      \diagram{d}{2.5em}{2.5em}{
        A_1 & A_2 & A_3 & A_4 & \cdots & A_{n - 1} & A_n & \mathbb{G}A_1\\
        A_1 & B_2 & B_3 & B_4 & \cdots & B_{n - 1} & B_n & \mathbb{G}
        A_1\\
        A_2 & B_2 & C_3 & C_4 & \cdots & C_{n - 1} & C_n & \mathbb{G}
        A_2\\
      };

      \path[->,midway,font=\scriptsize]
        (d-1-1) edge node[above] {$\underline{f_1}$} (d-1-2)
        ([xshift=-0.1em] d-1-1.south) edge[-] ([xshift=-0.1em] d-2-1.north)
        ([xshift=0.1em] d-1-1.south) edge[-] ([xshift=0.1em] d-2-1.north)
        (d-1-2) edge node[above] {$\underline{f_2}$} (d-1-3)
                     edge node[right] {$\underline{\varphi_2}$} (d-2-2)
        (d-1-3) edge node[above] {$\underline{f_3}$} (d-1-4)
                     %edge[densely dashed] node[right] {$\varphi_3$} (d-2-3)
        (d-1-4) edge node[above] {$\underline{f_4}$} (d-1-5)
                     %edge[densely dashed] node[right] {$\varphi_4$} (d-2-4)
        (d-1-5) edge node[above] {$\underline{f_{n - 2}}$} (d-1-6)
        (d-1-6) edge node[above] {$\underline{f_{n - 1}}$} (d-1-7)
                     %edge[densely dashed] node[right] {$\underline{\varphi_{n - 1}}$} (d-2-6)
        (d-1-7) edge node[above] {$\underline{\alpha_n}$} (d-1-8)
                     %edge[densely dashed] node[right] {$\varphi_n$} (d-2-7)
        ([xshift=-0.1em] d-1-8.south) edge[-] ([xshift=-0.1em] d-2-8.north)
        ([xshift=0.1em] d-1-8.south) edge[-] ([xshift=0.1em] d-2-8.north)
        (d-2-1) edge node[above] {$\underline{g_1}$} (d-2-2)
                     edge node[right] {$\underline{f_1}$} (d-3-1)
        (d-2-2) edge node[above] {$\underline{g_2}$} (d-2-3)
        ([xshift=-0.1em] d-2-2.south) edge[-] ([xshift=-0.1em] d-3-2.north)
        ([xshift=0.1em] d-2-2.south) edge[-] ([xshift=0.1em] d-3-2.north)
        (d-2-3) edge node[above] {$\underline{g_3}$} (d-2-4)
%                     edge[densely dashed] node[right] {$\underline{\theta_3}$} (d-3-3)
%        (d-2-4) edge node[above] {$\underline{g_4}$} (d-2-5)
%                     edge[densely dashed] node[right] {$\underline{\theta_4}$} (d-3-4)
        (d-2-5) edge node[above] {$\underline{g_{n - 2}}$} (d-2-6)
        (d-2-6) edge node[above] {$\underline{g_{n - 1}}$} (d-2-7)
%                     edge[densely dashed] node[right] {$\underline{\theta_{n - 1}}$} (d-3-6)
        (d-2-7) edge node[above] {$\underline{\beta_n}$} (d-2-8)
%                     edge[densely dashed] node[right] {$\underline{\theta_n}$} (d-3-7)
        (d-2-8) edge node[right] {$\mathbb{G}\underline{f_1}$} (d-3-8)
        (d-3-1) edge node[above] {$\underline{\varphi_2}$} (d-3-2)
        (d-3-2) edge node[above] {$\underline{h_2}$} (d-3-3)
        (d-3-3) edge node[above] {$\underline{h_3}$} (d-3-4)
        (d-3-4) edge node[above] {$\underline{h_4}$} (d-3-5)
        (d-3-5) edge node[above] {$\underline{h_{n - 2}}$} (d-3-6)
        (d-3-6) edge node[above] {$\underline{h_{n - 1}}$} (d-3-7)
        (d-3-7) edge node[above] {$\underline{\gamma_n}$} (d-3-8);
        %(d-1-4) edge[densely dashed,out=-102,in=30] node[pos=0.15,left] {$\underline{\psi_4}$}
%          (d-3-3)
%        (d-1-7) edge[densely dashed,out=-102,in=30] node[pos=0.15,left] {$\underline{\psi_n}$}
%          (d-3-6);
    \end{tikzpicture}
  \end{center}
  with commuting squares and with rows in $\Theta'$. Since $\underline{\varphi_{2}}\circ\underline{f_{1}}=\underline{g_{1}}$, there are two morphisms $u:A_{1}\rightarrow X$ and $v:X\rightarrow B_{2}$ such that $\varphi_{2} f_{1}-g_{1}=vu$. Since $g_{1}:A_{1}\rightarrow B_{2}$ is  an $\X$-monic, there exists a morphism
$w\colon B_2\to X$ such that $wg_{1}=u$. Thus, $\varphi_{2} f_{1}=(1+vw)g_{1}$. Note that $\underline{g_{1}}=\underline{(1+vw)g_{1}}$. Thus, there is a standard right $n$-angle
   $$A_1\xrightarrow{\underline{(1+vw)g_{1}}}B_2\xrightarrow{\underline{g_2}}B_3\xrightarrow{\underline{g_3}}\cdots\xrightarrow{\underline{g_{n-1}}}
B_n\xrightarrow{\underline{\beta_n}} \mathbb{G}A_1\in\Theta'.$$
By the construction of $\Theta'$, we have a right $n$-$\C$-sequence
$$A_1\xrightarrow{(1+vw)g_{1}}B_2\xrightarrow{g_2}B_3\xrightarrow{g_3}\cdots\xrightarrow{g_{n-1}} B_n\xrightarrow{g_n}  \Sigma A_{1}.$$
   By (RSN4), there exist the dotted
  morphisms  such that the following each square commutes
  \begin{center}
    \begin{tikzpicture}
      \diagram{d}{2.5em}{2.5em}{
        A_1 & A_2 & A_3 & A_4 & \cdots & A_{n - 1} & A_n & \Sigma A_1\\
        A_1 & B_2 & B_3 & B_4 & \cdots & B_{n - 1} & B_n & \Sigma
        A_1\\
        A_2 & B_2 & C_3 & C_4 & \cdots & C_{n - 1} & C_n & \Sigma
        A_2\\
      };

      \path[->,midway,font=\scriptsize]
        (d-1-1) edge node[above] {$f_1$} (d-1-2)
        ([xshift=-0.1em] d-1-1.south) edge[-] ([xshift=-0.1em] d-2-1.north)
        ([xshift=0.1em] d-1-1.south) edge[-] ([xshift=0.1em] d-2-1.north)
        (d-1-2) edge node[above] {$f_2$} (d-1-3)
                     edge node[right] {$\varphi_2$} (d-2-2)
        (d-1-3) edge node[above] {$f_3$} (d-1-4)
                     edge[densely dashed] node[right] {$\varphi_3$} (d-2-3)
        (d-1-4) edge node[above] {$f_4$} (d-1-5)
                     edge[densely dashed] node[right] {$\varphi_4$} (d-2-4)
        (d-1-5) edge node[above] {$f_{n - 2}$} (d-1-6)
        (d-1-6) edge node[above] {$f_{n - 1}$} (d-1-7)
                     edge[densely dashed] node[right] {$\varphi_{n - 1}$} (d-2-6)
        (d-1-7) edge node[above] {$f_n$} (d-1-8)
                     edge[densely dashed] node[right] {$\varphi_n$} (d-2-7)
        ([xshift=-0.1em] d-1-8.south) edge[-] ([xshift=-0.1em] d-2-8.north)
        ([xshift=0.1em] d-1-8.south) edge[-] ([xshift=0.1em] d-2-8.north)
        (d-2-1) edge node[above] {$(1+vw)g_{1}$} (d-2-2)
                     edge node[right] {$f_1$} (d-3-1)
        (d-2-2) edge node[above] {$g_2$} (d-2-3)
        ([xshift=-0.1em] d-2-2.south) edge[-] ([xshift=-0.1em] d-3-2.north)
        ([xshift=0.1em] d-2-2.south) edge[-] ([xshift=0.1em] d-3-2.north)
        (d-2-3) edge node[above] {$g_3$} (d-2-4)
                     edge[densely dashed] node[right] {$\theta_3$} (d-3-3)
        (d-2-4) edge node[above] {$g_4$} (d-2-5)
                     edge[densely dashed] node[right] {$\theta_4$} (d-3-4)
        (d-2-5) edge node[above] {$g_{n - 2}$} (d-2-6)
        (d-2-6) edge node[above] {$g_{n - 1}$} (d-2-7)
                     edge[densely dashed] node[right] {$\theta_{n - 1}$} (d-3-6)
        (d-2-7) edge node[above] {$g_n$} (d-2-8)
                     edge[densely dashed] node[right] {$\theta_n$} (d-3-7)
        (d-2-8) edge node[right] {$\Sigma f_1$} (d-3-8)
        (d-3-1) edge node[above] {$\varphi_2$} (d-3-2)
        (d-3-2) edge node[above] {$h_2$} (d-3-3)
        (d-3-3) edge node[above] {$h_3$} (d-3-4)
        (d-3-4) edge node[above] {$h_4$} (d-3-5)
        (d-3-5) edge node[above] {$h_{n - 2}$} (d-3-6)
        (d-3-6) edge node[above] {$h_{n - 1}$} (d-3-7)
        (d-3-7) edge node[above] {$h_n$} (d-3-8)
        (d-1-4) edge[densely dashed,out=-102,in=30] node[pos=0.15,left] {$\psi_4$}
          (d-3-3)
        (d-1-7) edge[densely dashed,out=-102,in=30] node[pos=0.15,left] {$\psi_n$}
          (d-3-6);
    \end{tikzpicture}
  \end{center}
 and the sequence
  \begin{center}
      \begin{tikzpicture}
        % Top row
        \node (A3) at (0,1.25){$A_3$};
        \node (A4B3) at (2,1.25){$A_4 \oplus B_3$};
        \node (A5B4C3) at (5,1.25){$A_5 \oplus B_4 \oplus C_3$};
        \node (A6B5C4) at (9,1.25){$A_6 \oplus B_5 \oplus C_4$};
        \node (tdots) at (12,1.25){$\cdots$};

        % Bottom row
        \node (mdots) at (0.25,0){$ $};
        \node (AnBn-1Cn-2) at (4.75,0){$A_n \oplus B_{n - 1} \oplus
          C_{n - 2}$};
        \node (BnCn-1) at (9.5,0){$B_n \oplus C_{n - 1}$};
        \node (mend) at (12.25,0){$C_n$};
        \node (mend2) at (14,0){$\Sigma A_3$};

        % Horizontal arrows
        \begin{scope}[font=\scriptsize,->,midway,above]
          % Top row
          \draw (A3) -- node{$\left[
              \begin{smallmatrix}
                f_3\\
                \varphi_3
              \end{smallmatrix}
            \right]$} (A4B3);
          \draw (A4B3) -- node{$\left[
              \begin{smallmatrix}
                -f_4 & 0\\
                \hfill \varphi_4 & -g_3\\
                \hfill \psi_4 & \hfill \theta_3
              \end{smallmatrix}
            \right]$} (A5B4C3);
          \draw (A5B4C3) -- node{$\left[
              \begin{smallmatrix}
                -f_5 & 0 & 0\\
                -\varphi_5 & -g_4 & 0\\
                \hfill \psi_5 & \hfill \theta_4 & h_3
              \end{smallmatrix}
            \right]$} (A6B5C4);
          \draw (A6B5C4) -- node{$\left[
              \begin{smallmatrix}
                -f_6 & 0 & 0\\
                \hfill \varphi_6 & -g_5 & 0\\
                \hfill \psi_6 & \hfill \theta_5 & h_4
              \end{smallmatrix}
            \right]$} (tdots);

          % Bottom row
          \draw (mdots) -- node{$\left[
              \begin{smallmatrix}
                -f_{n - 1} & 0 & 0\\
                (-1)^{n - 1} \varphi_{n - 1} & -g_{n - 2} & 0\\
                \psi_{n-1} & \theta_{n-2} & h_{n - 3}
              \end{smallmatrix}
            \right]$} (AnBn-1Cn-2);
          \draw (AnBn-1Cn-2) -- node{$\left[
              \begin{smallmatrix}
                (-1)^n \varphi_n & -g_{n - 1} & 0\\
                \psi_{n} & \theta_{n-1} & h_{n - 2}
              \end{smallmatrix}
            \right]$} (BnCn-1);
          \draw (BnCn-1) -- node{$\left[
              \begin{smallmatrix}
                \theta_{n} & h_{n - 1}
              \end{smallmatrix}
            \right]$} (mend);
          \draw (mend) -- node{$\Sigma f_2 \circ h_n$}
            (mend2);
        \end{scope}
      \end{tikzpicture}
    \end{center}
is a right $n$-$\C$-sequence in $R_n(\C,\Sigma)$ with $\left[
              \begin{smallmatrix}
                f_3\\
                \varphi_3
              \end{smallmatrix}
            \right]$ is an $\X$-monic.
Thus we have  the following commutative square diagram
 \begin{center}
    \begin{tikzpicture}
      \diagram{d}{2.5em}{2.5em}{
        A_1 & A_2 & A_3 & A_4 & \cdots & A_{n - 1} & A_n & \mathbb{G}A_1\\
        A_1 & B_2 & B_3 & B_4 & \cdots & B_{n - 1} & B_n & \mathbb{G}
        A_1\\
        A_2 & B_2 & C_3 & C_4 & \cdots & C_{n - 1} & C_n & \mathbb{G}
        A_2\\
      };

      \path[->,midway,font=\scriptsize]
        (d-1-1) edge node[above] {$\underline{f_1}$} (d-1-2)
        ([xshift=-0.1em] d-1-1.south) edge[-] ([xshift=-0.1em] d-2-1.north)
        ([xshift=0.1em] d-1-1.south) edge[-] ([xshift=0.1em] d-2-1.north)
        (d-1-2) edge node[above] {$\underline{f_2}$} (d-1-3)
                     edge node[right] {$\underline{\varphi_2}$} (d-2-2)
        (d-1-3) edge node[above] {$\underline{f_3}$} (d-1-4)
                     edge[densely dashed] node[right] {$\underline{\varphi_3}$} (d-2-3)
        (d-1-4) edge node[above] {$\underline{f_4}$} (d-1-5)
                     edge[densely dashed] node[right] {$\underline{\varphi_4}$} (d-2-4)
        (d-1-5) edge node[above] {$\underline{f_{n - 2}}$} (d-1-6)
        (d-1-6) edge node[above] {$\underline{f_{n - 1}}$} (d-1-7)
                     edge[densely dashed] node[right] {$\underline{\varphi_{n - 1}}$} (d-2-6)
        (d-1-7) edge node[above] {$\underline{\alpha_n}$} (d-1-8)
                     edge[densely dashed] node[right] {$\varphi_n$} (d-2-7)
        ([xshift=-0.1em] d-1-8.south) edge[-] ([xshift=-0.1em] d-2-8.north)
        ([xshift=0.1em] d-1-8.south) edge[-] ([xshift=0.1em] d-2-8.north)
        (d-2-1) edge node[above] {$\underline{g_1}$} (d-2-2)
                     edge node[right] {$\underline{f_1}$} (d-3-1)
        (d-2-2) edge node[above] {$\underline{g_2}$} (d-2-3)
        ([xshift=-0.1em] d-2-2.south) edge[-] ([xshift=-0.1em] d-3-2.north)
        ([xshift=0.1em] d-2-2.south) edge[-] ([xshift=0.1em] d-3-2.north)
        (d-2-3) edge node[above] {$\underline{g_3}$} (d-2-4)
                     edge[densely dashed] node[right] {$\underline{\theta_3}$} (d-3-3)
        (d-2-4) edge node[above] {$\underline{g_4}$} (d-2-5)
                     edge[densely dashed] node[right] {$\underline{\theta_4}$} (d-3-4)
        (d-2-5) edge node[above] {$\underline{g_{n - 2}}$} (d-2-6)
        (d-2-6) edge node[above] {$\underline{g_{n - 1}}$} (d-2-7)
                     edge[densely dashed] node[right] {$\underline{\theta_{n - 1}}$} (d-3-6)
        (d-2-7) edge node[above] {$\underline{\beta_n}$} (d-2-8)
                     edge[densely dashed] node[right] {$\underline{\theta_n}$} (d-3-7)
        (d-2-8) edge node[right] {$\mathbb{G}\underline{f_1}$} (d-3-8)
        (d-3-1) edge node[above] {$\underline{\varphi_2}$} (d-3-2)
        (d-3-2) edge node[above] {$\underline{h_2}$} (d-3-3)
        (d-3-3) edge node[above] {$\underline{h_3}$} (d-3-4)
        (d-3-4) edge node[above] {$\underline{h_4}$} (d-3-5)
        (d-3-5) edge node[above] {$\underline{h_{n - 2}}$} (d-3-6)
        (d-3-6) edge node[above] {$\underline{h_{n - 1}}$} (d-3-7)
        (d-3-7) edge node[above] {$\underline{\gamma_n}$} (d-3-8)
        (d-1-4) edge[densely dashed,out=-102,in=30] node[pos=0.15,left] {$\underline{\psi_4}$}
          (d-3-3)
        (d-1-7) edge[densely dashed,out=-102,in=30] node[pos=0.15,left] {$\underline{\psi_n}$}
          (d-3-6);
    \end{tikzpicture}
  \end{center}
  and the sequence
\begin{center}
      \begin{tikzpicture}
        % Top row
        \node (A3) at (0,1.25){$A_3$};
        \node (A4B3) at (2,1.25){$A_4 \oplus B_3$};
        \node (A5B4C3) at (5,1.25){$A_5 \oplus B_4 \oplus C_3$};
        \node (A6B5C4) at (9,1.25){$A_6 \oplus B_5 \oplus C_4$};
        \node (tdots) at (12,1.25){$\cdots$};

        % Bottom row
        \node (mdots) at (0.25,0){$ $};
        \node (AnBn-1Cn-2) at (4.75,0){$A_n \oplus B_{n - 1} \oplus
          C_{n - 2}$};
        \node (BnCn-1) at (9.5,0){$B_n \oplus C_{n - 1}$};
        \node (mend) at (12.25,0){$C_n$};
        \node (mend2) at (14,0){$\mathbb{G}A_3$};

        % Horizontal arrows
        \begin{scope}[font=\scriptsize,->,midway,above]
          % Top row
          \draw (A3) -- node{$\left[
              \begin{smallmatrix}
                \underline{f_3}\\
                \underline{\varphi_3}
              \end{smallmatrix}
            \right]$} (A4B3);
          \draw (A4B3) -- node{$\left[
              \begin{smallmatrix}
                \underline{-f_4} & 0\\
                \hfill \underline{\varphi_4} & \underline{-g_3}\\
                \hfill \underline{\psi_4} & \hfill \underline{\theta_3}
              \end{smallmatrix}
            \right]$} (A5B4C3);
          \draw (A5B4C3) -- node{$\left[
              \begin{smallmatrix}
                \underline{-f_5} & 0 & 0\\
                \underline{-\varphi_5} & \underline{-g_4} & 0\\
                \hfill \underline{\psi_5} & \hfill \underline{\theta_4} & \underline{h_3}
              \end{smallmatrix}
            \right]$} (A6B5C4);
          \draw (A6B5C4) -- node{$\left[
              \begin{smallmatrix}
                \underline{-f_6} & 0 & 0\\
                \hfill \underline{\varphi_6} & \underline{-g_5} & 0\\
                \hfill \underline{\psi_6 }& \hfill \underline{\theta_5} & \underline{h_4}
              \end{smallmatrix}
            \right]$} (tdots);

          % Bottom row
          \draw (mdots) -- node{$\left[
              \begin{smallmatrix}
                \underline{-f_{n - 1}} & 0 & 0\\
                \underline{(-1)^{n - 1} \varphi_{n - 1}} & -\underline{g_{n - 2}} & 0\\
                \underline{\psi_{n-1}} & \underline{\theta_{n-2}} & \underline{h_{n - 3}}
              \end{smallmatrix}
            \right]$} (AnBn-1Cn-2);
          \draw (AnBn-1Cn-2) -- node{$\left[
              \begin{smallmatrix}
                \underline{(-1)^n \varphi_n }& \underline{-g_{n - 1}} & 0\\
               \underline{ \psi_{n} }& \underline{\theta_{n-1}} & \underline{h_{n - 2}}
              \end{smallmatrix}
            \right]$} (BnCn-1);
          \draw (BnCn-1) -- node{$\left[
              \begin{smallmatrix}
                \underline{\theta_{n}} &\underline{ h_{n - 1}}
              \end{smallmatrix}
            \right]$} (mend);
          \draw (mend) -- node{$\underline{\lambda}$}
            (mend2);
        \end{scope}
      \end{tikzpicture}
    \end{center}
is a standard right $n$-angle in $\C/\X$. It suffices to show that $\underline{\lambda}=\mathbb{G}\underline{f_2}\circ \gamma_n$.

Consider the following commutative diagram
$$\xymatrix@C=1.8cm{A_2\ar[r]^{\varphi_2}\ar[d]^{f_2}&B_2\ar[r]^{h_2}\ar[d]^{\left[
              \begin{smallmatrix}
                0\\g_2
              \end{smallmatrix}
            \right]}&C_3\ar[r]^{h_3}\ar[d]^{\left[
              \begin{smallmatrix}
                0\\0\\1
              \end{smallmatrix}
            \right]}
&\cdot\cdot\cdot\ar[r]^{h_{n-1}}&C_n\ar[r]^{h_n}\ar@{=}[d]&\Sigma A_2\ar[d]^{\Sigma f_2}\\
A_3\ar[r]^{\left[
              \begin{smallmatrix}
                f_3\\ \varphi_3
              \end{smallmatrix}
            \right]}&A_4\oplus B_3\ar[r]^{\left[
              \begin{smallmatrix}
                -f_4 & 0\\
                \hfill \varphi_4 & -g_3\\
                \hfill \psi_4 & \hfill \theta_3
              \end{smallmatrix}
            \right]\quad}&A_5\oplus B_4\oplus C_3\ar[r]^{\quad\left[
              \begin{smallmatrix}
                -f_5 & 0 & 0\\
                -\varphi_5 & -g_4 & 0\\
                \hfill \psi_5 & \hfill \theta_4 & h_3
              \end{smallmatrix}
            \right]}&\cdot\cdot\cdot\ar[r]^{\left[
              \begin{smallmatrix}
                \theta_{n} & h_{n - 1}
              \end{smallmatrix}
            \right]}&C_n\ar[r]^{\Sigma f_2 \circ h_n}&\Sigma A_3}$$
of right $n$-$\C$-sequences in $R_n(\C,\Sigma)$, by Lemma \ref{y8}, we have
$\underline{\lambda}=\mathbb{G}\underline{f_2}\circ \gamma_n$.

(2) It is dual of (1).  \qed

\begin{remark}
In Theorem \ref{y9}, when $n=3$, it is just the Theorem 3.6 in \cite{[Li]}.
\end{remark}

\section{Frobenius $n$-Prile categories}
\setcounter{equation}{0}

In this section, we introduce the notion of a Frobenius $n$-prile category, and construct an $n$-angulated quotient category from an $n$-prile category.

\begin{definition}
Let $(\A,R_n(\C,\Sigma),\X)$ be a right $n$-suspended category and
 $(\A,L_n(\C,\Omega),\X)$  a  left $n$-suspended category. The quadruple
 $(\A,L_n(\C,\Omega),R_n(\C,\Sigma),\X)$ is called a \emph{ Frobenius $n$-prile category}
 if it satisfies the following conditions:
\begin{itemize}
\item[(a)] $(\Omega,\Sigma)$ is an adjoint pair with $\psi$ as the adjunction isomorphism.

\item[(b)] If
$\Omega A_n\xrightarrow{f_1}A_1\xrightarrow{f_2}A_2\xrightarrow{f_3}A_3\xrightarrow{f_4}\cdots
\xrightarrow{f_{n}}A_n$ is in $L_n(\C,\Omega)$ with $A_2,\cdots,A_n\in\C$, then
$$A_1\xrightarrow{f_2}A_2\xrightarrow{f_3}A_3\xrightarrow{f_4}\cdots\xrightarrow{f_{n}}A_n\xrightarrow{\psi_{A_n,A_1}(f_1)}\Sigma A_1$$ is in $R_{n}(\C,\Sigma)$. Conversely, if
$A_1\xrightarrow{g_1}A_2\xrightarrow{g_2}A_3\xrightarrow{g_3}\cdots\xrightarrow{g_{n-1}}A_n\xrightarrow{g_n}\Sigma A_1$ is in $R_{n}(\C,\Sigma)$ with $A_1,\cdots, A_{n-1}\in\C$, then
$$\Omega A_n\xrightarrow{\psi^{-1}_{A_n,A_1}(g_n)} A_1\xrightarrow{g_1}A_2\xrightarrow{g_2}A_3\xrightarrow{g_3}\cdots\xrightarrow{g_{n-1}}A_n$$ is in $L_{n}(\C,\Omega)$.
\item[(c)] For each object $A\in\C$, $x_{n-1}$ is a left $\X$-approximation of $A$ in the fixed right
$n$-$\C$-sequence $$A\xrightarrow{x}X_2\xrightarrow{x_2}X_3\xrightarrow{x_3}
\cdots\xrightarrow{x_{n-2}}X_{n-1}\xrightarrow{x_{n-1}}B\xrightarrow{y}\Sigma A,$$
where $X_2,\cdots, X_{n-1}\in\X$ and $x'_{1}$ is a right $\X$-approximation of $C$ in the fixed left $n$-$\C$-sequence $$\Omega A\xrightarrow{y'}C\xrightarrow{x'_{1}}X'_2\xrightarrow{x'_2}X'_3\xrightarrow{x'_3}
\cdots\xrightarrow{x'_{n-2}}X'_{n-1}\xrightarrow{x'_{n-1}}A,$$
where $X'_2,\cdots, X'_{n-1}\in\X$.
\end{itemize}
\end{definition}

\begin{example}\label{ex5.2}
\emph{(1)} Let $\B$ be a Frobenius $(n-2)$-exact category, and $\I$ the subcategory of projective-injective
objects of $\B$. Then $(\B,L_n(\B,0),R_{n}(\B,0),\I)$ is a Frobenius $n$-prile category, where
$L_n(\B,\Omega)$ and $R_{n}(\B,\Sigma)$ are constructed in Proposition \ref{p3}.

\emph{(2)} Let $\A$ be an $n$-angulated category and $\X$ two subcategories of $\C$.
If $\C$ is extension closed and $(\C,\C)$ is an $\X$-mutation pair, then
$(\A,L_n(\C,\Sigma^{-1}),R_{n}(\C,\Sigma),\X)$ is a Frobenius $n$-prile category, where
$L_n(\C,\Sigma^{-1})$ and $R_{n}(\C,\Sigma)$ are constructed in Example \ref{e100}.
\end{example}

\begin{lemma}\label{q2}
Let $(\A,L_n(\C,\Omega),R_n(\C,\Sigma),\X)$ be a Frobenius $n$-prile category.
Each solid commutative diagram
$$\xymatrix@C=0.7cm{
A\ar[r]^{a}\ar@{-->}[d]^{f} & X_2 \ar[r]^{x_2}\ar@{-->}[d]^{\varphi_2} & X_3 \ar[r]^{x_3}\ar@{-->}[d]^{\varphi_3} & \cdots\ar[r]^{x_{n-3}}&X_{n-2}\ar@{-->}[d]^{\varphi_{n-2}} \ar[r]^{x_{n-2}}&X_{n-1}\ar[r]^{x_{n-1}}\ar[d]^{^{\varphi_{n-1}}}& B \ar[r]^{b}\ar[d]^{g} & \Sigma A \ar@{-->}[d]^{\Sigma f}\\
A'\ar[r]^{a'} & X'_2 \ar[r]^{x'_2} & X'_3 \ar[r]^{x'_3} & \cdots \ar[r]^{x'_{n-3}}&X'_{n-2}\ar[r]^{x'_{n-2}}&X'_{n-1}\ar[r]^{x'_{n-1}} & B' \ar[r]^{b'}& \Sigma A'\\
}$$
with rows are right $n$-$\C$-sequences can be completed to a morphism of
right $n$-$\C$-sequences.
\end{lemma}

\proof Consider the following commutative diagram
$$\xymatrix@C=0.7cm{
A\ar[r]^{a} & X_2 \ar[r]^{x_2} & X_3 \ar[r]^{x_3}& \cdots\ar[r]^{x_{n-3}}&X_{n-2} \ar[r]^{x_{n-2}}&X_{n-1}\ar[r]^{x_{n-1}}\ar[d]^{^{\varphi_{n-1}}}& B \ar[r]^{b}\ar[d]^{g} & \Sigma A \\
A'\ar[r]^{a'} & X'_2 \ar[r]^{x'_2} & X'_3 \ar[r]^{x'_3} & \cdots \ar[r]^{x'_{n-3}}&X'_{n-2}\ar[r]^{x'_{n-2}}&X'_{n-1}\ar[r]^{x'_{n-1}} & B' \ar[r]^{b'}& \Sigma A'\\
}$$
of right $n$-$\C$-sequences. By assumption, we have the following commutative
$$\xymatrix@C=1.2cm{
\Omega B\ar[r]^{\psi^{-1}_{B,A}(b)}\ar[d]^{\Omega(g)}&A\ar[r]^{a}\ar@{-->}[d]^{f} & X_2 \ar[r]^{x_2}\ar@{-->}[d]^{\varphi_2} & X_3 \ar[r]^{x_3}\ar@{-->}[d]^{\varphi_3} & \cdots\ar[r]^{x_{n-3}}&X_{n-2}\ar@{-->}[d]^{\varphi_{n-2}} \ar[r]^{x_{n-2}}&X_{n-1}\ar[r]^{x_{n-1}}\ar[d]^{^{\varphi_{n-1}}}& B \ar[d]^{g}\\
\Omega B'\ar[r]^{\psi^{-1}_{B',A'}(b')}&A'\ar[r]^{a'} & X'_2 \ar[r]^{x'_2} & X'_3 \ar[r]^{x'_3} & \cdots \ar[r]^{x'_{n-3}}&X'_{n-2}\ar[r]^{x'_{n-2}}&X'_{n-1}\ar[r]^{x'_{n-1}} & B'
}$$
of left $n$-$\C$-sequences, where $f,\varphi_2,\cdots,\varphi_{n-2}$
exist by the axiom (LSN3) of a  left $n$-suspended category.

We claim that $\Sigma f\circ b=b'g$. Indeed, by the naturality of $\psi^{-1}_{B',A'}$ in $B'$,
we have the following commutative diagram
$$\xymatrix@C=1.2cm{
\Hom_{\A}(B',\Sigma A')\ar[r]^{\psi^{-1}_{B',A'}}\ar[d]^{g^{\ast}}
&\Hom_{\A}(\Omega B',A')\ar[d]^{\Omega(g)^{\ast}} \\
\Hom_{\A}(B,\Sigma A')\ar[r]^{\psi^{-1}_{B,A'}}&\Hom_{\A}(\Omega B,A')
.}$$
It follows that for the morphism $b'\colon B'\to \Sigma A'$, we have
$$\psi^{-1}_{B',A'}(b)\circ\Omega(g)=\psi^{-1}_{B,A'}(b'g).$$
In a similar way,
by the naturality of $\psi^{-1}_{B,A}$ in $B$,
we have the following commutative diagram
$$\xymatrix@C=1.2cm{
\Hom_{\A}(B,\Sigma A)\ar[r]^{\psi^{-1}_{B,A}}\ar[d]^{(\Sigma f)_{\ast}}
&\Hom_{\A}(\Omega B,A)\ar[d]^{f_{\ast}} \\
\Hom_{\A}(B,\Sigma A')\ar[r]^{\psi^{-1}_{B,A'}}&\Hom_{\A}(\Omega B,A')
.}$$
It follows that for the morphism $b\colon B\to \Sigma A$, we have
$$\psi^{-1}_{B,A'}(\Sigma f\circ b)=f\circ\psi^{-1}_{B,A}(b).$$
Since $$\psi^{-1}_{B',A'}(b')\circ\Omega(g)=f\circ\psi^{-1}_{B,A}(b),$$
we have $$\psi^{-1}_{B,A'}(b'g)=\psi^{-1}_{B,A'}(\Sigma f\circ b)$$and then
$\Sigma f\circ b=b'g$. Thus we have the following commutative diagram
$$\xymatrix@C=0.7cm{
A\ar[r]^{a}\ar@{-->}[d]^{f} & X_2 \ar[r]^{x_2}\ar@{-->}[d]^{\varphi_2} & X_3 \ar[r]^{x_3}\ar@{-->}[d]^{\varphi_3} & \cdots\ar[r]^{x_{n-3}}&X_{n-2}\ar@{-->}[d]^{\varphi_{n-2}} \ar[r]^{x_{n-2}}&X_{n-1}\ar[r]^{x_{n-1}}\ar[d]^{^{\varphi_{n-1}}}& B \ar[r]^{b}\ar[d]^{g} & \Sigma A \ar@{-->}[d]^{\Sigma f}\\
A'\ar[r]^{a'} & X'_2 \ar[r]^{x'_2} & X'_3 \ar[r]^{x'_3} & \cdots \ar[r]^{x'_{n-3}}&X'_{n-2}\ar[r]^{x'_{n-2}}&X'_{n-1}\ar[r]^{x'_{n-1}} & B' \ar[r]^{b'}& \Sigma A'\\
}$$
with rows are right $n$-$\C$-sequences.
\qed

\begin{theorem}\label{h11}
Let $(\A,L_n(\C,\Omega),R_n(\C,\Sigma),\X)$ be a Frobenius $n$-prile category. Then
the quotient category $\C/\X$ is an $n$-angulated category.
\end{theorem}

\proof \proof (1) By Theorem \ref{y9}, we know that $(\C/\X,\mathbb{G},\Theta')$
is a right $n$-angulated category. It suffices to show that $\mathbb{G}$ is an equivalence functor.

$\bullet$ Let us prove that $\mathbb{G}$ is dense.

For any object $A\in\C/\X$, there exists a right $n$-$\C$-sequence
$$A\xrightarrow{a}X_2\xrightarrow{x_2}X_3\xrightarrow{x_3}
\cdots\xrightarrow{x_{n-2}}X_{n-1}\xrightarrow{x_{n-1}}C\xrightarrow{b}\Sigma C,$$
where $a$ is a left $\X$-approximation of $A$ and $X_2,X_3,\cdots, X_{n-1}\in\X$.
By the definition of $\mathbb{G}$ and Lemma \ref{y7}, we have
$C\simeq\mathbb{G}(A)$. This shows that $\mathbb{G}$ is dense.

$\bullet$ Let us prove that $\mathbb{G}$ is full.

For any morphism $\underline{h}\colon \mathbb{G}A\to \mathbb{G}B$ in $\C/\X$, since $(\A,L_n(\C,\Omega),R_n(\C,\Sigma),\X)$ is a Frobenius $n$-prile category, there are
right $n$-$\C$-sequences
$$\xymatrix@C=0.7cm{
A\ar[r]^{f_1}& X_2 \ar[r]^{f_2} & X_3 \ar[r]^{f_3}& \cdots \ar[r]^{f_{n-2}}&X_{n-1}\ar[r]^{f_{n-1}}& \mathbb{G}A \ar[r]^{f_n} & \Sigma A}$$
where $f_{n-1}$ is a right $\X$-approximation and $X_2,X_3,\cdots,X_{n-1}\in\X$.
$$\xymatrix@C=0.7cm{B \ar[r]^{g_1} & X'_2 \ar[r]^{g_2} & X'_3 \ar[r]^{g_3} & \cdots \ar[r]^{g_{n-2}}&X'_{n-1}\ar[r]^{g_{n-1}} & \mathbb{G}B \ar[r]^{g_n}& \Sigma B
}$$
where $g_{n-1}$ is a right $\X$-approximation and $X'_2,X'_3,\cdots,X'_{n-1}\in\X$.
Thus there exists the following commutative diagram of right $n$-$\C$-sequences
$$\xymatrix@C=0.7cm{
A \ar[r]^{f_1}\ar[d]^{f} & X_2 \ar[r]^{f_2}\ar[d]^{\varphi_2} & X_3 \ar[r]^{f_3}\ar[d]^{\varphi_3} & \cdots \ar[r]^{f_{n-2}}&X_{n-1}\ar[r]^{f_{n-1}}\ar[d]^{^{\varphi_{n-1}}}& \mathbb{G}A \ar[r]^{f_n}\ar[d]^{h} & \Sigma A \ar[d]^{\Sigma f}\\
B\ar[r]^{g_1} & X'_2 \ar[r]^{g_2} & X'_3 \ar[r]^{g_3} & \cdots \ar[r]^{g_{n-2}}&X'_{n-1}\ar[r]^{g_{n-1}} & \mathbb{G}B \ar[r]^{g_n}& \Sigma B\\
}$$
where $\varphi_{n-1}$ exists since $g_{n-1}$ is a right $\X$-approximation,
$f,\varphi_{1},\cdots,\varphi_{n-2}$ exist by Lemma \ref{q2}.
It follows that $\underline{h}=\mathbb{G}\underline{f}$. This shows that $\mathbb{G}$ is full.

$\bullet$ Let us prove that $\mathbb{G}$ is faithful.

Let $f\colon A\to B$ be a morphism in $\C$.  There exists the following commutative diagram
of right $n$-$\C$-sequences
$$\xymatrix@C=0.7cm{
A \ar[r]^{f_1}\ar[d]^{f} & X_2 \ar[r]^{f_2}\ar[d]^{\varphi_2} & X_3 \ar[r]^{f_3}\ar[d]^{\varphi_3} & \cdots \ar[r]^{f_{n-2}}&X_{n-1}\ar[r]^{f_{n-1}}\ar[d]^{^{\varphi_{n-1}}}& \mathbb{G}A \ar[r]^{f_n}\ar[d]^{h} & \Sigma A \ar[d]^{\Sigma f}\\
B \ar[r]^{g_1} & X'_2 \ar[r]^{g_2} & X'_3 \ar[r]^{g_3} & \cdots \ar[r]^{g_{n-2}}&X'_{n-1}\ar[r]^{g_{n-1}} & \mathbb{G}B \ar[r]^{g_n}& \Sigma B\\
}$$
where $f_{n-1}$ and $g_{n-1}$ are right $\X$-approximations.
Assume that $\mathbb{G}\underline{f}=\underline{h}=0$.
Then $h$ factors through some object in $\X$.
Since $g_{n-1}$ is a right $\X$-approximation of $\mathbb{G}B_1$. we have that
$h$ factors through $g_{n-1}$. Since $(\A,L_n(\C,\Omega),\X)$ be a left $n$-suspended category, we know that $f$ factors through $f_1$. Therefore $\underline{f}=0$. This shows that $\mathbb{G}$ is faithful.  \qed

\begin{remark}
In Theorem \ref{h11}, when $n=3$, it is just the Theorem 5.5 in \cite{[Li]}. Moreover, under  the situation  of Example \ref{ex5.2}, one can reobtain \cite[Theorem 5.11]{[J]} and \cite[Theorem 3.7]{[L2]} by Theorem \ref{h11}.
\end{remark}

\end{document}